\documentclass{article}
\usepackage{amsmath,amssymb,amscd,amsthm,amsfonts,amstext,amsxtra,
mathrsfs
}
\usepackage{eufrak,fontenc,fontsmpl,latexsym,rawfonts}
\theoremstyle{plain}
\newtheorem{theorem}{Theorem}
\newtheorem{lemma}[theorem]{Lemma}
\newtheorem{proposition}[theorem]{Proposition}
\newtheorem{corollary}[theorem]{Corollary}
\newtheorem{definition}[theorem]{Definition}

\newtheorem*{theorem*}{Theorem}
\newtheorem*{lemma*}{Lemma}
\newtheorem*{proposition*}{Proposition}
\newtheorem*{corollary*}{Corollary}
\newtheorem*{definition*}{Definition}
\newtheorem*{conjecture*}{Conjecture}
\theoremstyle{definition}
\newtheorem*{remark*}{Remark}

\newtheorem*{keywords}{Keywords}

\begin{document}
\title{Direct image for multiplicative and relative
$K$-theories from transgression of the families index theorem, part 2}
\author{Alain Berthomieu\\\and Institut de Math\'ematiques de Toulouse, UMR 5219\and Centre Universitaire de Formation et de Recherche \and Jean-Fran\c{c}ois Champollion, Campus d'Albi,\and Place de Verdun,
81012 Albi Cedex, France.
\and{\tt{alain.berthomieu@univ-jfc.fr}}}
\maketitle
\begin{abstract}
This is the sequel of \cite{MoiPartie1}, (the second part of the ``longstanding forthcoming preprint'' referred to in \cite{MoiOberwolfach}).

Here, the procedure of transgressing the families index
theorem (the so-called $\eta$-form) is adapted to take in account the case of Dirac type
operators with kernels of varying dimension.

The constructed form is then used to define the direct image under proper
submersions of the ``free multiplicative'' $K$-theory which
was defined in \cite{MoiPartie1}, the behaviour of the characteristic classes 
on free multiplicative $K$-theory defined in \cite{MoiPartie1}
under submersion is studied.

Some universal
caracterisation of the forms is provided.

Finally, combining our result with Bismut
and Lott's results
on direct images of flat
vector bundles yields a ``Grothendieck-Riemann-Roch''
theorem for Nadel-Chern-Simons classes on relative $K$-theory for
flat vector bundles which was defined in \cite{MoiPartie1}
(and whose
holomorphic counterpart is known from the previous work \cite{BerthomieuPNCDIRKT}).
\begin{keywords} multiplicative $K$-theory, families index, Chern-Simons transgression, proper submersions.
\end{keywords}
\noindent{\textbf{AMS-classification:}} Primary: 14F05, 19E20, 57R20, secondary: 14F40, 19D55, 53C05, 55R50.
\end{abstract}
\section{Introduction:}
In \cite{MoiPartie1} \S2.3,
a ``free multiplicative'' $K$-theory $\widehat K_{\rm{ch}}(M)$ inspired by \cite{KaroubiHCKT},
\cite{KaroubiTGCCS}, \cite{KaroubiCCFFHA}
and suggested by U. Bunke \cite{BunkeSchick}
was constructed for any smooth manifold $M$, whose objects are triples $(E,\nabla,\alpha)$ where $E$ is a complex vector
bundle on $M$ with connection $\nabla$ and $\alpha$ is an odd-degree differential form on $M$
defined modulo exact forms, with following relations: direct sum, and
change of connection in the following way:
$\ (E,\nabla_{\!0},\alpha)=\big(E,\nabla_{\!1},\alpha+\widetilde{\rm{ch}}(\nabla_{\!0},\nabla_{\!1})\big)$, where
$\widetilde{\rm{ch}}(\nabla_{\!0},\nabla_{\!1})$ is Chern ans Simons'
transgression form: it is defined modulo exact forms and verifies
\begin{equation}\label{Sarah1}d\widetilde{\rm{ch}}(\nabla_{\!0},\nabla_{\!1})={\rm{ch}}
(\nabla_{\!1})-{\rm{ch}}(\nabla_{\!0})\end{equation}
where ${\rm{ch}}
(\nabla_{\!0})$ and ${\rm{ch}}(\nabla_{\!1})$ are Chern-Weil representatives of the Chern character of $E$ computed from the connections $\nabla_{\!0}$ and $\nabla_{\!1}$.

The Chern character of $(E,\nabla,\alpha)$ is defined to be $\widehat{\rm{ch}}(E,\nabla,\alpha)=
{\rm{ch}}(\nabla)-d\alpha$
(it is a differential form, not only a cohomology class).
A Borel-type class ${\mathfrak{B}}(E,\nabla,\alpha)$ is also defined to be
$\widetilde{\rm{ch}}(\nabla^*,\nabla)-2i\, {\mathfrak{Im}}\, \alpha$
(this is a purely imaginary odd differential form defined modulo exact forms)
where $\nabla^*$ is the adjoint transpose of the connection
$\nabla$ with respect to any hermitian metric on $E$ (see \cite{MoiPartie1} \S2.3.2 and \S2.3.3).

This free multiplicative
$K$-theory with differential form Chern character is a theory which
together provides classes for vector bundles with connections and has
direct images under proper submersions (see corollary \ref{Cor:main} infra).

This Chern-Simons type transgression is generalised to transgress the family
index theorem in the following sense: 
let $\pi\colon M\to B$ be a proper smooth submersion between compact manifolds,
there is a (classical) direct image morphism $\pi_*^{\rm{Eu}}\colon K^0_{\rm{top}}(M)\longrightarrow K^0_{\rm{top}}(B)$ associated to the fibral Euler-Dirac $d+d^*$ operator, which is extensively studied in \cite{MoiPartie1}
\S3.1. In particular, for any vector bundle $\xi$ on $M$, there are constructions of representatives of $\pi_*^{\rm{Eu}}([\xi])$ of the form $[E^+]
-[E^-]\in K^0_{\rm{top}}(B)$ which are ``linked'' in the sense that there is a canonical isotopy class of vector bundle isomorphism
$E^+\oplus G^-\oplus H\overset\sim\longrightarrow E^-\oplus G^+\oplus H$ if $G^+$ and $G^-$ are another couple of constructed representatives (and $H$ is some vector bundle on $B$). 
Let $TZ$ be the vector bundle on $M$
consisting of vertical tangent vectors (those lying in the kernel of $d\pi$),
let $\nabla_{\!TZ}$ be some Levi-Civita connexion on $TZ$ and $e(\nabla_{\!TZ})$
be the Chern-Weil representative of the Euler class of $TZ$, the following result is proved in \S\S\ref{IranzoSarah} and \ref{S.I.} below:
\begin{theorem*}
To any vector bundles $\xi$ on $M$ with connection
$\nabla_{\!\xi}$ and $F^+$ and $F^-$ on $B$ with respective connections $\nabla_{\!F^+}$
and $\nabla_{\!F^-}$ such that $[F^+]-[F^-]=\pi_*^{\rm{Eu}}[\xi]$ in $K^0_{\text{top}}(B)$ with
some link between $F^+-F^-$ and any explicit representative of $\pi_*[\xi]$
obtained by the analytic families index construction of \cite{MoiPartie1} \S3.1, there exists a way to associate
a differential form $\tau$ on $B$ such that
\[d\tau={\rm{ch}}(\nabla_{\!F^+})-{\rm{ch}}(\nabla_{\!F^-})-\int_{M/B}e(\nabla_{\!TZ}){\rm{ch}}
(\nabla_{\!\xi})\]
This construction is functorial by pullbacks on fibered products, the obtained
form is real if $\nabla_{\!\xi}$, $\nabla_{\!F^+}$ and $\nabla_{\!F^-}$ all respect some
hermitian metrics on $\xi$, $F^+$ and $F^-$ respectively. It is aditive under direct sums, and it vanishes if ${\rm{dim}}M-{\rm{dim}}B$ is odd and $F^+$ and $F^-$ are taken to be null vector bundles with canonical link with constructed representatives from
\cite{MoiPartie1} \S3.2.3 lemma 10.
\end{theorem*}
$\int_{M/B}$ is integration on the fibres of $\pi$. The fact that the right hand side of the preceding equation is an exact form
is a consequence of the cohomological counterpart of the families index theorem
(see \cite{AtiyahSinger} theorem (5.1))

This construction is the classical construction of $\eta$-forms as in \cite{BismutCheeger} (see also \cite{BerlineGetzlerVergne} \S10.5), and the first goal of this paper is to
make it compatible with the procedures of constructing representatives of $\pi^{\rm{Eu}}_*[\xi]$ when the dimensions of the kernels of vertical Dirac operators vary, and with the links (see \cite{MoiPartie1}\S3.1).
The notation $\tau$ initially chosen for this $\eta$-form is kept here
to avoid telescoping with other data called $\eta$. Of course this $\tau$ is not a torsion form (as in \cite{BismutLott} definition 3.22).

Afterwards, it is possible to exploit its various properties.
For example 
\begin{corollary}\label{Cor:main}
The map $(\xi,\nabla_{\!\xi},\alpha)\longmapsto\left(F^+\!,\!\nabla_{\!F^+},\int_{M/B}e(\nabla_{\!TZ})\alpha\right)
-(F^-\!,\!\nabla_{\!F^-},\tau)$ provides a morphism $\pi^{\rm{Eu}}_!\colon\widehat K_{\rm{ch}}(M)\to\widehat K_{\rm{ch}}(B)$ which
vanishes if ${\rm{dim}}M-{\rm{dim}}B$ is odd and is such that

$\widehat{\rm{ch}}\big(\pi^{\rm{Eu}}_!(\xi,\nabla_{\!\xi},\alpha)\big)=\int_{M/B}e(\nabla_{\!TZ})\widehat{\rm{ch}}(\xi
\nabla_{\!\xi},\alpha)$

${\mathfrak{B}}\big(\pi^{\rm{Eu}}_!(\xi,\nabla_{\!\xi},\alpha)\big)=\int_{M/B}e(\nabla_{\!TZ})
{\mathfrak{B}}(\xi,\nabla_{\!\xi},\alpha)$
\end{corollary}
This last statement needs the reality property above.
The functoriality property allows to establish anomaly formulas for $\tau$
from which the dependence of $\pi^{\rm{Eu}}_!$ on metric datas (in $\nabla_{\!TZ}$
and in the construction of $\tau$) is explained: let $\widetilde e(\nabla_{\!TZ,0}
,\nabla_{\!TZ,1})$ denote
the Chern-Simons form associated to the Euler form $e$ and two connections
$\nabla_{\!TZ,0}$ and $\nabla_{\!TZ,1}$ on $TZ$, then the difference between
direct images of $(\xi,\nabla_{\!\xi},\alpha)$ constructed with metric datas
corresponding to connections $\nabla_{\!TZ,0}$ and $\nabla_{\!TZ,1}$ on
$TZ$ equals $\left(0,0,\int_{M/B}\widetilde e(\nabla_{\!TZ,0},\nabla_{\!TZ,1})
\widehat{\rm{ch}}(\xi,\nabla_{\!\xi},\alpha)\right)$. This in particular proves the independence
on metric datas of the restriction of the morphism $\pi^{\rm{Eu}}_!$ to the subgroup $MK^0(M)$
of $\widehat K_{\rm{ch}}(M)$
consisting of objects whose Chern character is a constant integer. This subgroup
is a Karoubi multiplicative $K$-theory group \cite{KaroubiHCKT}, \cite{KaroubiTGCCS}, \cite{KaroubiCCFFHA} (see lemma \ref{l:anomaliemetrique} infra). 

Consider a flat complex vector bundle $\xi$ on $M$ (with flat connection $\nabla_{\!\xi}$)
the even/odd degree de Rham cohomology of the fibres of $\pi$ with coefficients
in $\xi$ provide flat complex vector bundles $H^\pm$ on $B$ with flat connexions
$\nabla_{\!\pm}$. This provides a direct image $K^0_{\text{flat}}(M)\to
K^0_{\text{flat}}(B)$ where $K^0_{\text{flat}}$ is the algebraic $K$-theory of
flat vector bundles modulo exact sequences. This morphism is compatible with
forgetful map $K^0_{\text{flat}}\to K^0_{\text{top}}$ and direct image for
topological $K$-theory. It is proved in lemma \ref{lemmederealite} below (see also \cite{MaZhang}) that
$\tau$ vanishes for data $(\xi,\nabla_{\!\xi})$, $(H^+,\nabla_{\!+})$ and $(H^-,\nabla_{\!-})$,
from which one obtains the following
\begin{corollary*}
The direct images on $K^0_{\text{flat}}$ and on $\widehat K_{\rm{ch}}$ are compatible
with the morphism $(\xi,\nabla_{\!\xi})\in K^0_{\text{flat}}\longmapsto(\xi,\nabla_{\!\xi},0)
\in\widehat K_{\rm{ch}}$.
%
\end{corollary*}

In \cite{MoiPartie1} \S2.1.4, a relative $K$-theory group $K^0_{\rm{rel}}(M)$ was constructed
whose objects are quintuples $(E,\nabla_{\!E},F,\nabla_{\!F},f)$ where
$E$ and $F$ are vector bundles with flat connections $\nabla_{\!E}$
and $\nabla_{\!F}$ and $f\colon E\overset\sim\longrightarrow
F$ is a smooth isomorphism, and a direct image map $\pi_*\colon
K^0_{\rm{rel}}(M)\longrightarrow K^0_{\rm{rel}}(B)$ was defined.
It was also proved that the assignment $(E,\nabla_{\!E},F,\nabla_{\!F},f)\longmapsto\widetilde{\rm{ch}}
(\nabla_{\!E},f^*\nabla_{\!F})$ provides a characteristic class morphism
${\mathcal N}_{\rm{ch}}\colon K^0_{\rm{rel}}(M)\longmapsto H^{\rm{odd}}
(M,{\mathbb C})$.

The vanishing property of $\tau$ for flat bundles, its functoriality
and its anomaly formulas
allow to prove the following (see \S\ref{GRRR} below)
\begin{corollary}\label{CorollSarah} for any $(E,F,f)\in K^0_{\text{rel}}(M)$:
\[{\mathcal N}_{\rm{ch}}\big(\pi_*(E,\nabla_{\!E},F,\nabla_{\!F},f)\big)=\int_{M/B}e(\nabla_{\!TZ}){\mathcal N}_{\rm{ch}}(E,\nabla_{\!E},F,\nabla_{\!F},f)\]
(This is a cohomological formula).
\end{corollary}

The paper is organized as follows: in section \ref{eta}, the construction
of $\tau$ is performed, its various properties are studied,
and it is characterized when possible, it is then applied in section \ref{cinq} to prove the above corollaries. I plan to deal with double submersions in a third (and last) part. Of course constant reference will be made to \cite{MoiPartie1} in the sequel.
\section{Transgression of the families index theorem}\label{eta}
\subsection{Bismut and Lott's ``Levi-Civita'' superconnection:}\label{s:geomdata}
The fibres of $\pi$ are supposed to be modelled on the compact manifold $Z$. Let $\xi$ be a complex vector bundle on $M$ with connection $\nabla_{\!\xi}$ and hermitian metric $h^\xi$.
The infinite rank vector bundle ${\mathcal E}^\pm=C^\infty(Z,\wedge^{\genfrac{}{}{0pt}{}{\rm{even}}{\rm{odd}}}T^*Z\otimes\xi)$ (of fibral $\xi$-valued differential forms) on $B$ was defined in \cite{MoiPartie1} \S3.1.1. If $TZ$ is endowed with some
euclidean metric $g^Z$, then $\wedge^\bullet T^*Z\otimes\xi$ is endowed with a hermitian metric $(\ /\ )_Z$ and
${\mathcal E}$ with a $L^2$ scalar product $\langle\ \, ,\ \rangle_{L^2}$ as defined in
\cite{MoiPartie1} \S3.1.1.

Let $T^H\!M$ denote some supplementary subvector bundle to $TZ$ in $TM$
(as in \cite{MoiPartie1} \S3.3.3), for any vector ${\tt u}$ tangent to $B$ at $y$, its horizontal lift ${\tt u}^H$ on the fibre $\pi^{-1}(y)$
is the section of $T^H\!M$ over $\pi^{-1}(y)$ such that $d\pi({\tt u}^H)={\tt u}$ at any point of $\pi^{-1}(y)$.
Call $P^{TZ}$ the projection of $TM$ onto $TZ$ with kernel
$T^H\!M$. If ${\tt u}$ and ${\tt v}$ are vector fields defined on a neighbourhood of
$y\in B$, then the vector field $P^{TZ}[{\tt u}^H,{\tt v}^H]$ on $Z_y=\pi^{-1}(y)$ depends
on the values of ${\tt u}$ and ${\tt v}$ at $y$ only.
Let $\iota_T\colon\wedge^2TB\longrightarrow{\text{End}}^{\text{odd}}
({\mathcal E})$ be the operator
which to ${\tt u}$ and ${\tt v}\in T_yB$ associates the interior product by
$-P^{TZ}[{\tt u}^H,{\tt v}^H]$ in $\wedge^\bullet T^*Z\otimes\xi$.
$\iota_T$ can be extended to a globally odd ${\rm{End}}{\mathcal E}$-valued differential form (of differential form degree $2$) on $B$.

Let $d^{\nabla_{\!\xi}}$ be the exterior differential operator associated to $\nabla_{\!\xi}$ along the fibres of $\pi$, it is considered here as an endomorphism of ${\mathcal E}$.
The connection $\overline\nabla$ on the bundle ${\mathcal E}$
was defined in \cite{MoiPartie1} \S3.3.3. Then
$\overline\nabla+d^{\nabla_{\!\xi}}+\iota_T$ is a superconnection on
${\mathcal E}$
in the sense of Quillen \cite{Quil}, \cite{BerlineGetzlerVergne} Definitions 1.37 and 9.12
and \cite{Bismut}.
It can be proved to coincide with the total exterior differential operator
$d^{M}$ on $\xi$-valued differential forms (defined using $\nabla_{\!\xi}$) on $M$ (throw the identification
of formula (37) in \cite{MoiPartie1}) as in \cite{BerlineGetzlerVergne} proposition 10.1
(the proof of \cite{BismutLott} \S III (b) cannot be adapted here because
$(d^{M})^2$ does not vanish if $\nabla_{\!\xi}$ is not flat).

Let $T\wedge\colon\wedge^2TB\longrightarrow{\text{End}}^{\text{odd}}
({\mathcal E})$ be the operator which
associates to ${\tt u}$ and ${\tt v}\in T_yB$ the exterior product in 
${\mathcal E}_y$ by the one form $(-P^{TZ}[{\tt u}^H,{\tt v}^H])^\flat$
(the dual throw $g^{Z}$ to
the vector field $-P^{TZ}[{\tt u}^H,{\tt v}^H]$) on $\pi^{-1}(y)$.
Before and after being extended to a globally odd ${\rm{End}}{\mathcal E}$-valued differential form on $B$ (of differential form degree $2$), $T\wedge$ is the adjoint of $\iota_T$, so that $\iota_T-T\wedge$ is a special autoadjoint ${\rm{End}}{\mathcal E}$-valued differential form in the sense of \cite{MoiPartie1} \S2.2.4.

$d^{\nabla_{\!\xi}}$ and its adjoint $(d^{\nabla_{\!\xi}})^*$ with respect to $\langle\ \, ,\ \rangle_{L^2}$ are also mutually special adjoint (in the sense of \cite{MoiPartie1} \S2.2.4) as
${\rm{End}}{\mathcal E}$-valued differential forms (with differential form degree $0$). The adjoint connection $\overline\nabla^S$ to $\overline\nabla$ with respect to $\langle\ \, ,\ \rangle_{L^2}$ was defined in \cite{MoiPartie1} \S3.3.3.
The superconnection
$\overline\nabla^S+(d^{\nabla_{\!\xi}})^*
-T\wedge$ is the adjoint transpose of the superconnection
$\overline\nabla
+d^{\nabla_{\!\xi}}+\iota_T$ in the sense of \cite{BismutLott} \S I(b) and
proposition 3.7, (and also of \cite{MoiPartie1} formula (28)).

The relevant Bismut-Levi-Civita superconnection in this
context is defined for any $t>0$ as in \cite{BismutLott} (3.50) (and also (3.49),
(3.30) and proposition 3.4):
\begin{equation}\label{LCSK}
C_t=\frac12\big(\overline\nabla+\overline\nabla^S\big)+\frac{\sqrt{t}}2\big(d^{\nabla_{\!\xi}}+(d^{\nabla_{\!\xi}})^*\big)+
\frac{1}{2\sqrt{t}}(\iota_T-T\wedge)\end{equation}
In the case of a fibered product of $\widetilde B$ and $M$ over $B$ as
in \cite{MoiPartie1} formula (39):
\begin{equation}\label{pullbackfunctorial}
\begin{CD}\widetilde B\times_B\!M@>>>M\\@VVV@VVV\\\widetilde B@>>>B\end{CD}
\end{equation}
(the model ofthe fibre may not change),
the construction of $C_t$ is functorial if
the horizontal subspace $T^H\!(\widetilde B\times_B\!M)$ is taken to be the
subspace of $T(\widetilde B\times_B\!M)$ consisting of vectors which are sent by the tangent map of $\widetilde B\times_B\!M\longrightarrow M$ into $T^H\!M$. (It is not always the pullback of $T^H\!M$). 
\subsection{Properties of its Chern character:}
$C_t^2$ is a fiberwise positive second order elliptic
differential operator so that its heat kernel $\exp-C_t^2$
is trace class. The Chern character of $C_t$ is defined to be
\[{\text{ch}}(C_t)=\phi{\text{Tr}}_s\exp-C_t^2\]
where ${\rm{Tr}}_s$ is a supertrace, namely the trace on ${\mathcal E}^+$ minus the trace on ${\mathcal E}^-$, the ${\mathbb Z}_2$-graduation being given by the parity of the degree of vertical differential forms, and $\phi$ divides $2k$-degree horizontal differential forms by $(2\pi i)^k$ (see \cite{MoiPartie1} \S2.2.4).
Let $\nabla^*_{\!\xi}$ be the adjoint connection to $\nabla_{\!\xi}$
with respect to $h^\xi$ (in the sense of \cite{BismutLott}, see
formula (9) in \cite{MoiPartie1})
\begin{lemma}\label{pascoole!}
${\rm{ch}}(C_t)$ is a real form. Let $C_t\!\!\check{\ }$ be the corresponding superconnection on ${\mathcal E}$
constructed from $\nabla_{\!\xi}^*$ instead of $\nabla_{\!\xi}$
then:
\[{\rm{ch}}(C_t\!\!\check{\ })=(-1)^{{\rm{dim}}Z}{\rm{ch}}(C_t)\]
In particular, ${\rm{ch}}(C_t)$ vanishes if ${\rm{dim}}Z$ is odd and $\nabla_{\!\xi}$ respects the hermitian metric $h^\xi$. Moreover ${\rm{ch}}(C_t)$ is a constant integer if $\nabla_{\!\xi}$ is flat.
\end{lemma}
\begin{proof}
The reality of ${\rm{ch}}(C_t)$ is due to the fact that $C_t$ is for any $t$
the half sum of $\overline\nabla+\sqrt{ t}d^{\nabla_{\!\xi}}+\frac1{\sqrt t}\iota_T$ and $\overline\nabla^S+\sqrt{t}(d^{\nabla_{\!\xi}})^*-\frac1{\sqrt t}T\wedge$ which are mutually adjoint transpose superconnections, so that the
considerations of the end of \S2.2.4 of \cite{MoiPartie1} apply here.

The case of flat connection $\nabla_{\!\xi}$ is treated in \cite{BismutLott} Theorem 3.15.

For any ${\tt w}\in TZ$, its Clifford action $c({\tt w})$ on $\wedge^\bullet T^*Z$ and the Hodge-Clifford operator $*_Z$ 
(see \cite{MoiPartie1} \S3.2.1 for the definition of both)
commute if ${\rm{dim}}Z$ is odd and
anticommute if ${\rm{dim}}Z$ is even. Thus $\iota_T-T\wedge=-c(T)$
behaves similarly.
It follows then from formula (45) in \cite{MoiPartie1} that \begin{align*}
\frac{\sqrt{t}}2\big(d^{\nabla_{\!\xi}}+(d^{\nabla_{\!\xi}})^*\big)+
\frac{1}{2\sqrt{t}}&(\iota_T-T\wedge)=\\=-(-1)^{{\rm{dim}}Z}*_Z^{-1}&\left(
\frac{\sqrt{t}}2\big(d^{\nabla^*_{\!\xi}}+(d^{\nabla^*_{\!\xi}})^*\big)+
\frac{1}{2\sqrt{t}}(\iota_T-T\wedge)\right)*_Z
\end{align*}

The formulas (55) of \cite{MoiPartie1} and their obvious counterpart
for $\overline\nabla\check{\ }$ and $\overline\nabla\check{\ }^S$
prove that $\frac12(\overline\nabla\check{\ }+\overline\nabla\check{\ }^S)
=\frac12(\overline\nabla+\overline\nabla^S)$. This connection
will be called $\overline\nabla^u$ in the sequel.

It is a consequence of \cite{MoiPartie1} lemma 12 that the covariant derivative with respect
to $\overline\nabla^u$ commutes with $*_Z$. For ${\mathbb Z}_2$-graduation reasons,
this proves that the exterior derivative associated with $\overline\nabla^u$
on ${\rm{End}}E$-valued differential forms on $B$ supercommutes with $*_Z$.

Let $N_H$ be the graduation operator on $\wedge T^*B$ which multiplies $k$-degree
differential forms by $k$, the properties above give the following formulas:
\begin{equation}\label{people}
\begin{aligned}
C_t&=-(-1)^{{\rm{dim}}Z}(-1)^{N_H}*_Z^{-1}C_t\!\!\check{\ }*_Z(-1)^{N_H}\\
{\text{so that }}\qquad\qquad C_t^2&=(-1)^{N_H}*_Z^{-1}C_t\!\!\check{\ }^2*_Z(-1)^{N_H}\\
{\text{and }}\qquad \exp(-C_t^2)&=(-1)^{N_H}*_Z^{-1}\exp(-C_t\!\!\check{\ }^2)*_Z(-1)^{N_H}
\end{aligned}
\end{equation}
But $\exp(-C^2_t)$ and $\exp(-C_t\!\!\check{\ }^2)$ are globally even
${\rm{End}}{\mathcal E}$-valued differential form, so that only their even
differential form degree parts contribute to its supertrace.

In this context of infinite rank vector bundles, it remains true that the supertrace
of the supercommutator of two $L^2$-bounded ${\rm{End}}{\mathcal E}$-valued differential forms,
one of which is trace class, vanishes. Applying this to $[*_Z^{-1},\omega*_Z]$
gives the following relation
\begin{equation}\label{starsupertrace}
{\rm{Tr}}_s(\omega)=(-1)^{{\rm{dim}}Z}{\rm{Tr}}_s(*_Z^{-1}\omega*_Z)
\end{equation}
valid for any globally even ${\rm{End}}{\mathcal E}$-valued trace-class differential form $\omega$.

The equation of the lemma is then a consequence of this relation, commutation
rules \eqref{people} and the comment after it.
\end{proof}

Put $\nabla_{\!\xi}^u=\frac12(\nabla_{\!\xi}+\nabla_{\!\xi}^*)$, then
\begin{proposition}\label{locindthm}As $t$ tends to $0$, ${\rm{ch}}(C_t)$
has for any $k\geq1$ an asymptotic of the form
\[{\rm{ch}}(C_t)=\sum_{j=0}^{k-1}t^{j+\frac12}A_j+{\mathcal O}(t^{k+\frac12})\qquad{\text{ if }}{\rm{dim}}Z{\text{ is odd}}\]
\[{\rm{ch}}(C_t)=\sum_{j=0}^{k-1}t^{j}B_j+{\mathcal O}(t^{k})\qquad{\text{ if }}{\rm{dim}}Z{\text{ is even}}\]
in either case:
\[\mathop{\text{lim}}\limits_{t\to 0}{\rm{ch}}(C_t)
=\int_Ze(\nabla_{\!TZ})\wedge{\rm{ch}}
(\nabla^u_{\!\xi})\]
\end{proposition}
\begin{proof}The asymptotics with $\sum\limits_{j=-\frac12{\rm{dim}}Z}^{k-1}$
are classical results on heat kernels (see \cite{BerlineGetzlerVergne} \S\S2.5 and 2.6
and appendix after \S9.7).

The limit formula (and thus the vanishing of the terms $A_j$ and/or $B_j$
for negative $j$) is a
consequence of \cite{BismutLott} (3.76). The connection $\nabla_{\!\xi}$
is supposed to be flat in \cite{BismutLott}, which is not the case here: thus formula
\cite{BismutLott} (3.52) does not hold true here.

However, consider ${\mathscr R}$
defined as
in \cite{BismutLott} (3.56) without taking \cite{BismutLott} (3.52) into account,
then
the $z=0$ case of the Lichnerowicz-type formula
of \cite{BismutLott} theorem 3.11 holds true here.

Thus the rescaling formula \cite{BismutLott} (3.75) and its consequence \cite{BismutLott}
(3.76) remain true here.
(This is only a matter of Clifford degrees which has nothing to do with
the fact that $\nabla_{\!\xi}$ be flat or not).

In particular, if ${\rm{dim}}Z$ is odd, then
the same argument as in \cite{BismutLott} (3.79) applies,
and both sides of the equality vanish.
\end{proof}
\subsection{Adapting the Bismut-Levi-Civita superconnection
to the varying dimension of the kernel of $d^{\nabla_{\!\xi}}+(d^{\nabla_{\!\xi}})^*$:}
\label{s:adapt}
Consider some smooth real increasing function $\chi$ on
${\mathbb R}_+$ which vanishes on $[0,\frac12]$ and equals $1$ on $[1,+\infty)$.
Consider some suitable data $(\eta^+,\eta^-,\psi)$ with respect to $\pi$ in the sense of \S3.1 of \cite{MoiPartie1}.
Put some hermitian metrics $h^\pm$ on $\eta^\pm$ and some connection
$\nabla_{\!\eta}$ on $\eta^+\oplus\eta^-$ which respects the decomposition.
Consider the following $t$-depending
superconnection on $({\mathcal E}^+\oplus\eta^+)\oplus({\mathcal E}^-\oplus\eta^-)$:
\begin{equation}\label{BLLCmodif}B_t=\overline\nabla^u\oplus\nabla_{\!\eta}
+\frac{\sqrt{t}}2{\mathcal D}^{\nabla_{\!\xi}}_{\chi(t)\psi}
+\frac{1}{2\sqrt{t}}(\iota_T-T\wedge)=C_t\oplus\nabla_{\!\eta}+
\frac{\sqrt{t}}2\chi(t)
\left(\psi+\psi^*\right)\end{equation}
(the modified fiberwise Euler-Dirac operator ${\mathcal D}^{\nabla_{\!\xi}}_{\chi(t)\psi}$ is defined as in \cite{MoiPartie1} \S3.1.1). $B_t^2$ is as $C^2_t$ a fiberwise positive
second order elliptic operator, so that its heat kernel is trace class.
Its Chern character is defined in the same way, the supertrace to consider being the trace on ${\text{End}}({\mathcal E}^+\oplus\eta^+)$
minus the trace on ${\text{End}}({\mathcal E}^-\oplus\eta^-)$.
For $t\leq\frac12$, one has
${\rm{ch}}(B_t)={\rm{ch}}(C_t)+{\rm{ch}}(\nabla_{\!\eta})$.

$\psi$ is of differential form degree $0$ so that $\psi^*$ is the special adjoint of $\psi$. It follows from an argument of the same kind as was used in lemma \ref{pascoole!} that ${\rm{ch}}(B_t^2)$ is for any $t$ a real form provided $\nabla_{\!\eta}$
respects $h^+$ and $h^-$.

Consider the orthogonal projection $P^{{\mathcal H}^\pm}\colon{\mathcal E}^\pm
\oplus\eta^\pm\longrightarrow {\mathcal H}^\pm$ on  ${\mathcal H}^\pm={\rm{Ker}}{\mathcal D}^\pm_\psi$, (and $P^{\mathcal H}=P^{{\mathcal H}^+}
\oplus P^{{\mathcal H}^-}$). The associated connection
on ${\mathcal H}={\mathcal H}^+\oplus{\mathcal H}^-$ is
\begin{equation}\label{nablahache}\nabla_{\!\mathcal H}=P^{\mathcal H}\left(
\overline\nabla^u\oplus
\nabla_{\!\eta}\right)P^{\mathcal H}\end{equation}

This connection respects the decomposition ${\mathcal H}^+\oplus{\mathcal H}^-$, and it also respects the hermitian metric on ${\mathcal H}$ obtained
by restriction of $\langle\ \, ,\, \ \rangle_{L^2}\oplus h^\pm$ provided
$\nabla_{\!\eta}$ respects $h^\pm$ (as can be proved by a direct elementary computation).

It is proved in \cite{BerlineGetzlerVergne} theorem 9.26 that:
\begin{equation}\label{limiteinfinie}\mathop{\text{lim}}\limits_{t\to+\infty}{\text{ch}}(B_t)=
{\text{ch}}(\nabla_{\!\mathcal H})\end{equation}
in the sense of any ${\mathcal C}^\ell$ norm on any compact subset of
$B$.

\begin{lemma}\label{LovelySarah}
If ${\rm{dim}}Z$ is even, let $B_t\!\!\check{\ }$ be the as above modified superconnection constructed
from $C_t\!\!\check{\ }$ (or $\nabla_{\!\xi}^*$) and the suitable data
of \cite{MoiPartie1} \S3.2.2, namely $\big(\eta^+,\eta^-,-(*_Z\oplus{\rm{Id}}_{\eta^-})\circ\psi\circ
(*_Z^{-1}\oplus{\rm{Id}}_{\eta^+})\big)$. Then
\[\phi
{\text{Tr}}_s(\exp-B_t\!\!\check{\ }^2)=\phi
{\text{Tr}}_s(\exp-B_t^2)\]

If ${\rm{dim}}Z$ is odd, let $B_t\!\!\check{\ }$ be the as above modified superconnection constructed
from $C_t\!\!\check{\ }$ (or $\nabla_{\!\xi}^*$) and the suitable data
$\big(\eta^-, \eta^+,(*_Z\oplus{\rm{Id}}_{\eta^+})\circ
\psi^*\circ(*_Z^{-1}\oplus{\rm{Id}}_{\eta^-})\big)$
(as at the end of \cite{MoiPartie1} \S3.2.3). Then
\[\phi
{\text{Tr}}_s(\exp-B_t\!\!\check{\ }^2)=-\phi
{\text{Tr}}_s(\exp-B_t^2)\]
\end{lemma}
\begin{proof}
If ${\rm{dim}}Z$ is even, denote $-(*_Z\oplus{\rm{Id}}_{\eta^-})\circ\psi\circ
(*_Z^{-1}\oplus{\rm{Id}}_{\eta^+})$ by $\Psi$, the suitable data $(\eta^+,\eta^-,\Psi)$ for $\nabla_{\!\xi}^*$
are chosen so that formula (45) of \cite{MoiPartie1} reads now
\begin{align*}
{\mathcal D}^{\nabla_{\!\xi}+}_\psi&=-(*_Z\oplus{\rm{Id}}_{\eta^-})^{-1}
{\mathcal D}^{\nabla^*_{\!\xi}+}_\Psi(*_Z\oplus{\rm{Id}}_{\eta^+})\\
{\mathcal D}^{\nabla_{\!\xi}-}_\psi&=-(*_Z\oplus{\rm{Id}}_{\eta^+})^{-1}
{\mathcal D}^{\nabla^*_{\!\xi}-}_\Psi(*_Z\oplus{\rm{Id}}_{\eta^-})
\end{align*}

If ${\rm{dim}}Z$ is odd, denote $(*_Z\oplus{\rm{Id}}_{\eta^+})\circ
\psi^*\circ(*_Z^{-1}\oplus{\rm{Id}}_{\eta^-})$ by $\Psi$, the suitable data $(\eta^-,\eta^+,\Psi)$ for $\nabla_{\!\xi}^*$
are chosen so that formula (45) of \cite{MoiPartie1} reads now
\begin{align*}
{\mathcal D}^{\nabla_{\!\xi}+}_\psi&=(*_Z\oplus{\rm{Id}}_{\eta^-})^{-1}
{\mathcal D}^{\nabla^*_{\!\xi}-}_\Psi(*_Z\oplus{\rm{Id}}_{\eta^+})\\
{\mathcal D}^{\nabla_{\!\xi}-}_\psi&=(*_Z\oplus{\rm{Id}}_{\eta^+})^{-1}
{\mathcal D}^{\nabla^*_{\!\xi}+}_\Psi(*_Z\oplus{\rm{Id}}_{\eta^-})
\end{align*}
Put ${\rm{Id}}_\eta={\rm{Id}}_{\eta^+}\oplus{\rm{Id}}_{\eta^-}$, this
has the following consequences for $B_t$ (for any $t$):
\begin{equation}\label{HighClassSarah}
B_t=-(-1)^{{\rm{dim}}Z}(-1)^{N_H}(*_Z\oplus{\rm{Id}}_\eta)^{-1}B_t\!\!\check{\ }(*_Z\oplus{\rm{Id}}_\eta)(-1)^{N_H}\end{equation}
The end of the proof is similar to the end of the proof of lemma \ref{pascoole!}
(from the first line of \eqref{people} on, which is here replaced by this last equation).
\end{proof}
Both $B_t$ and its Chern
character are functorial by pullbacks on fibered products as in
\eqref{pullbackfunctorial} (if the horizontal subspace of the
source manifold is taken as described at the end of \S\ref{s:geomdata}).

Note also that the construction can be performed with any smooth function $\chi$
on $B\times{\mathbb R}_+$ which vanishes on $B\times[0,\varepsilon]$ and equals
$1$ on $B\times[A,+\infty)$ for any $0<\varepsilon<A$, and which is increasing
with respect to the variable in ${\mathbb R}_+$. This is of course not
essential, but will be useful to prove some independence of the constructed
forms on the choice of the function $\chi$.
\subsection{Properties of the infinitesimal transgression form:}\label{s:4.4}
Consider the following classical transgression formula (see \cite{MoiPartie1} (19))
\[\frac{d}{dt}{\text{ch}}(B_t)=-d\left[\phi
{\text{Tr}}_s\left(\frac{\partial B_t}{\partial t}
\exp-B^2_t\right)\right]\]
so that for any $0<S<T<+\infty$
\begin{equation}\label{transgressionadistancefinie}{\text{ch}}(B_S)-
{\text{ch}}(B_T)=d\left[\int_S^T\phi
{\text{Tr}}_s\left(\frac{\partial B_t}{\partial t}
\exp-B^2_t\right)dt\right]\end{equation}
\begin{lemma}\label{l:reel}$\phi{\rm{Tr}}_s\left(\frac{\partial B_t}{\partial t}
\exp-B^2_t\right)$ is a real form if $\nabla_{\!\eta}$ respects $h^\pm$ (the hermitian metrics on $\eta^\pm$, if not, this form is changed into its complex conjugate if $\nabla_{\!\eta}$ is changed into its adjoint transpose connection
with respect to $h^\pm$).

Let $B_t\!\!\check{\ }$ be as in lemma \ref{LovelySarah}. Then
\[\phi
{\text{Tr}}_s\left(\frac{\partial B_t\!\!\check{\ }}{\partial t}
\exp-B_t\!\!\check{\ }^2\right)=(-1)^{{\rm{dim}}Z}\phi
{\text{Tr}}_s\left(\frac{\partial B_t}{\partial t}
\exp-B_t^2\right)\]
In particular, $\phi
{\text{Tr}}_s\left(\frac{\partial C_t}{\partial t}
\exp-C_t^2\right)$ vanishes if $\nabla_{\!\xi}$ respects $h^\xi$ and ${\rm{dim}}Z$ is odd.
\end{lemma}
\begin{proof}
In one hand, $\exp-B^2_t$ is a globally even 
${\rm{End}}{\mathcal E}$-valued differential form, so that its supercommutator with $\frac{\partial B_t}{\partial t}$
is their usual commutator; and it is special autoadjoint in the sense of \cite{MoiPartie1} \S2.2.4 if
$\nabla_{\!\eta}$ respects $h^\pm$ on $\eta^\pm$ (if not, the two forms
obtained from mutually transpose adjoint connections on $\eta$ are mutually special adjoint).

In the other hand, $\frac{\partial B_t}{\partial t}$ is for any $t$ a special autoadjoint ${\rm{End}}{\mathcal E}$-valued differential form
in the sense of \cite{MoiPartie1} \S2.2.4, so that the product $\frac{\partial B_t}{\partial t}
\exp-B^2_t$ is the special adjoint of $(\exp-B^2_t)\frac{\partial B_t}{\partial t}$ (if $\nabla_{\!\eta}$ respect $h^\pm$). Thus
\[\phi{\rm{Tr}}_s\left(\frac{\partial B_t}{\partial t}
\exp-B^2_t\right)=\phi{\rm{Tr}}_s\left((\exp-B^2_t)\frac{\partial B_t}{\partial t}\right)=\overline{\phi{\rm{Tr}}_s\left(\frac{\partial B_t}{\partial t}
\exp-B^2_t\right)}\]
and the reality follows (the case when $\nabla_{\!\eta}$ does not respect
$h^\pm$ is similar).

As was made in \eqref{people}, it successively follows from \eqref{HighClassSarah} that
\begin{align*}
\exp(-B_t^2)&=(-1)^{N_H}(*_Z\oplus{\rm{Id}}_\eta)^{-1}\exp(-B_t\!\!\check{\ }^2)(*_Z\oplus{\rm{Id}}_\eta)(-1)^{N_H}\\
\frac{\partial B_t}{\partial t}\exp(-B_t^2)&=\\=-(-1&)^{{\rm{dim}}Z}(-1)^{N_H}(*_Z\oplus{\rm{Id}}_\eta)^{-1}\frac{\partial B_t\!\!\check{\ }}{\partial t}\exp(-B_t\!\!\check{\ }^2)(*_Z\oplus{\rm{Id}}_\eta)(-1)^{N_H}
\end{align*}
$\frac{\partial B_t}{\partial t}\exp(-B_t^2)$ and $\frac{\partial B_t\!\!\check{\ }}{\partial t}\exp(-B_t\!\!\check{\ }^2)$ are globally odd, so that only their odd differential form degree parts contribute to their
supertrace. And because in both cases ${\rm{Id}}_\eta$ has the same
parity as $*_Z$, the counterpart of \eqref{starsupertrace} for a globally odd
${\rm{End}}{\mathcal E}$-valued differential form $\omega$ reads here
\[{\rm{Tr}}_s(\omega)={\rm{Tr}}_s\big((*_Z\oplus{\rm{Id}}_\eta)^{-1}\omega(*_Z\oplus{\rm{Id}}_\eta)\big)
\] 
The last assertions of the lemma follow.
\end{proof}
\subsection{Asymptotics of the infinitesimal transgression form:}

Consider now the product manifold $\widetilde M={\mathbb R}\times M$ and its obvious submersion
$\widetilde\pi={\rm{Id}}_{\mathbb R}\times\pi$ onto $\widetilde B={\mathbb R}\times B$. Extend $\xi$ tautologically to $\widetilde M$
with constant (with respect to $s$) hermitian metric and connection $d_{\mathbb R}+\nabla_{\!\xi}$ (where $d_{\mathbb R}=ds\frac\partial{\partial s}$ is the trivial canonical differential along ${\mathbb R}$).
Consider any smooth real positive function $f$ on ${\mathbb R}$ such that $f(1)=1$ and
endow the vertical tangent bundle of $\widetilde\pi$ with the metric
$\frac1{f(s)}g^Z$. Choose $T^H\!\widetilde M=T{\mathbb R}\oplus T^H\!M$
as horizontal bundle of $\widetilde\pi$. Let's calculate the
Bismut-Lott Levi-Civita superconnection $\widetilde C_t$ in this context.

The equivalent here of the connection $\overline\nabla$ defined in
\cite{MoiPartie1} (53) is obviously equal to $d_{\mathbb R}+\overline\nabla$.
The vertical exterior differential operator $d^{\nabla_{\!\xi}}$ is unchanged, and so is the operator $\iota_T$ (defined at the beginning of \S\ref{s:geomdata}).

The volume form of the fibres of $\widetilde\pi$ on $\{s\}\times B$ is equal to
$f(s)^{-\frac{\rm{dimZ}}2}$ times the corresponding volume form of the fibres on $\{1\}\times B$. The ponctual scalar product between vertical differential forms
of degree $k$ on $\{s\}\times B$ is equal to the one on $\{1\}\times B$
multiplied by $f(s)^k$. 
Call $\widetilde{\mathcal E}$
the infinite rank vector bundle on $\widetilde B$ of $\xi$-valued vertical
differential forms, and define $N_V\in{\rm{End}}{\mathcal E}$ or ${\rm{End}}\widetilde{\mathcal E}$ to be the operator which multiplies vertical differential forms
by their degree. The global $L^2$ scalar product on the restriction of
$\widetilde{\mathcal E}$ to $\{s\}\times B$ is thus equal to $f(s)^{N_V-\frac{{\rm{dim}}Z}2}\langle\ \, ,\ \rangle_{L^2}$ (where
$\langle\ \, ,\ \rangle_{L^2}$ defined in \cite{MoiPartie1} (38) is the one on
$\{1\}\times B$).

It follows that the adjoint of $d^{\nabla_{\!\xi}}$ is $f(s)(d^{\nabla_{\!\xi}})^*$ (if $(d^{\nabla_{\!\xi}})^*$ is its adjoint on $\{1\}\times B$) and the adjoint of $\iota_T$ is $\frac1{f(s)}T\wedge$ (if $T\wedge$ is its adjoint on $\{1\}\times B$). In the same way, as in \cite{MoiPartie1} (54), one has $(d_{\mathbb R}+\overline\nabla)^S
=d_{\mathbb R}+ds\frac{f'(s)}{f(s)}\big(N_V-\frac{{\rm{dim}}Z}2\big)+\overline\nabla^S$.

Thus if $\widetilde C_{t,s}$ denotes the Bismut-Lott Levi-Civita superconnection on $\{s\}\times B$:
\begin{align*}
C_{t,s}&=\frac12\big(\overline\nabla+\overline\nabla^S\big)+\frac{\sqrt{t}}2\Big(d^{\nabla_{\!\xi}}+f(s)(d^{\nabla_{\!\xi}})^*\Big)+
\frac{1}{2\sqrt{t}}\left(\iota_T-\frac1{f(s)}T\wedge\right)\\
\widetilde C_t\, &=\, C_{t,s}+d_{\mathbb R}+\frac12ds\frac{f'(s)}{f(s)}\left(N_V-\frac{{\rm{dim}}Z}2\right)
\end{align*}

One then computes:
\begin{align*}
[d_{\mathbb R},C_{t,s}]&=\frac{\sqrt{t}}2f'(s)ds(d^{\nabla_{\!\xi}})^*
+\frac{f'(s)}{2\sqrt{t}f(s)^2}ds\, T\wedge\\
[N_V,C_{t,s}]&=\frac{\sqrt{t}}2
\big(d^{\nabla_{\!\xi}}-f(s)(d^{\nabla_{\!\xi}})^*\big)+\frac1{2\sqrt{t}}\left(-\iota_T-\frac1{f(s)}T\wedge\right)\\
\left[d_{\mathbb R}+\frac12ds\frac{f'(s)}{f(s)}N_V\, ,\, C_{t,s}\right]&=
\frac{\sqrt{t}}4ds\frac{f'(s)}{f(s)}\big(d^{\nabla_{\!\xi}}+f(s)(d^{\nabla_{\!\xi}})^*\big)+\\
&\qquad\qquad+ds\frac{f'(s)}{4\sqrt{t}f(s)}\left(-\iota_T+\frac1{f(s)}T\wedge\right)\\
&=t\,  ds\frac{f'(s)}{f(s)}\frac{\partial C_{t,s}}{\partial t}\end{align*}
\begin{align*}
\widetilde C_t^2=C_{t,s}^2+&\left[d_{\mathbb R}+\frac12ds\frac{f'(s)}{f(s)}\left(N_V-\frac{{\rm{dim}}Z}2\right)\, ,\, C_{t,s}\right]\\=C_{t,s}^2+&t\, ds\frac{f'(s)}{f(s)}\frac{\partial C_{t,s}}{\partial t}\end{align*}
\begin{equation}\label{mere}
{\rm{Tr}}_s\exp(-\widetilde C^2_t)={\rm{Tr}}_s\exp(-C_{t,s}^2)-
t\, ds\frac{f'(s)}{f(s)}{\rm{Tr}}_s\left(\frac{\partial C_{t,s}}{\partial t}
\exp(-C_{t,s}^2)\right)
\end{equation}
\begin{proposition*}If $\nabla_{\!\xi}$ is flat and if the suitable data used in the construction of $B_t$ are the trivial ones $(\{0\},\{0\},0)$ (see \cite{MoiPartie1} \S3.3.2), then:
\[\phi
{\text{Tr}}_s\left(\frac{\partial B_t}{\partial t}
\exp-B^2_t\right)=\phi
{\text{Tr}}_s\left(\frac{\partial C_t}{\partial t}
\exp-C^2_t\right)=0\]
In general, one has the following estimates
\begin{align*}{\rm{as }}\quad\, t\to+\infty:\qquad\qquad&\phi
{\text{Tr}}_s\left(\frac{\partial B_t}{\partial t}
\exp-B^2_t\right)={\mathcal O}(t^{-\frac32})\\{\rm{as }}\quad\, t\to0:\qquad\qquad&\phi
{\text{Tr}}_s\left(\frac{\partial B_t}{\partial t}
\exp-B^2_t\right)=\left\{\begin{aligned}{\mathcal O}&(1)\quad{\text{ if ${\rm{dim}}Z$ is even}}\\{\mathcal O}&(t^{-\frac12})\, {\text{ if ${\rm{dim}}Z$ is odd}}\end{aligned}\right.
\end{align*}
\end{proposition*}
\begin{proof}The first assertion is proved in \cite{MaZhang}.
It is reproved here as a direct consequence of \eqref{mere}, of the last assertion of lemma \ref{pascoole!} (and the fact that if $\nabla_{\!\xi}$ is
flat on $\xi$ over $M$, then $d_{\mathbb R}+\nabla_{\!\xi}$ is also flat on the pullback of $\xi$ over $\widetilde M$).

The $t\to+\infty$ asymptotic is proved by the adaptation of \cite{BerlineGetzlerVergne} theorem 9.23
which is proposed (though not detailed) at the end of \S9.3 of \cite{BerlineGetzlerVergne}. (Here $\chi(t)$ is constant on a neighbourhood of $+\infty$, so that the arguments of the proof of theorems 9.7 and 9.23 of \cite{BerlineGetzlerVergne} apply).

The second one will be proved with the technique proposed in \cite{BerlineGetzlerVergne} Theorem 10.32: apply proposition \ref{locindthm}
on $\widetilde M$, one obtains because of the factor $t$ appearing in \eqref{mere} an asymptotic of the form
\[\phi
{\text{Tr}}_s\left(\frac{\partial B_t}{\partial t}
\exp-B^2_t\right)=\left\{\begin{aligned}\sum_{j=-1}^{k-1}&E_jt^j+{\mathcal O}(t^k)\quad{\text{ if ${\rm{dim}}Z$ is even}}\\\sum_{j=0}^{k-1}&E_jt^{j-\frac12}+{\mathcal O}(t^{k-\frac12})\quad{\text{ if ${\rm{dim}}Z$ is odd}}\end{aligned}\right.\]
This proves the assertion for odd dimensional fibres.

Suppose now that ${\rm{dim}}Z$ is even. Let $\widetilde\nabla_{\!TZ}$ be the Levi-Civita connection on the vertical tangent bundle of the submersion $\widetilde\pi$ over $\widetilde M$
(It is the (unique) connection on this bundle obtained from the projection of $T\widetilde M$ on it with kernel $T^H\!\widetilde M$
and any Levi-Civita connection on $T\widetilde M$ corresponding to
any riemannian metric which coincides with $g^Z$ on vertical vectors
and make vertical vectors orthogonal to $T^H\!\widetilde M$, see
\cite{MoiPartie1}, beginning of the proof of lemma 12).
Denote by $\int_Z$ the integral along the fibres of $\widetilde\pi$, then $E_{-1}$ is the factor of $ds$ in the 
decomposition of the differential form
$\int_Ze(\widetilde\nabla_{\!TZ}){\rm{ch}}(\nabla_{\!\xi}^u)$ with respect to 
$\Omega(\widetilde B,{\mathbb C})=C^\infty\big({\mathbb R},\Omega(B,{\mathbb C})\big)\oplus ds\wedge C^\infty\big({\mathbb R},\Omega(B,{\mathbb C})\big)$. (Here and throughout, for any manifold $X$ and any vector bundle $\zeta$ on $X$, $\Omega(X,\zeta)$ is the space of $C^\infty$ $\zeta$-valued differential forms on $X$).
This is because the Chern character is functorial by pullbacks. However, $\widetilde\nabla_{\!TZ}$ is not the pullback of $\nabla_{\!TZ}$. A direct calculation from the classical formula for Levi-Civita connections (see
\cite{BerlineGetzlerVergne} formula (1.18)) yields
\[\widetilde\nabla_{\!TZ}=d_{\mathbb R}+\nabla_{\!TZ}+\frac{f'(s)}{2fs)}ds\]
so that $\widetilde\nabla_{\!TZ}^2=\nabla_{\!TZ}^2$ because $d_{\mathbb R}$ and $ds$ both commute with $\nabla_{\!TZ}$. Thus the curvature of $\widetilde\nabla_{\!TZ}$ is the pullback of the one of $\nabla_{\!TZ}$
and neither $e(\widetilde\nabla_{\!TZ})$ nor ${\rm{ch}}(\nabla_{\!\xi}^u)$
have a $ds$ component. This proves the vanishing of $E_{-1}$.
\end{proof}
These estimates together with formulae \eqref{limiteinfinie} and \eqref{transgressionadistancefinie} and proposition \ref{locindthm} provide the following transgression formula:
\[d\left[\int_0^{+\infty}\phi
{\text{Tr}}_s\left(\frac{\partial B_t}{\partial t}
\exp-B^2_t\right)dt\right] =\int_Ze(\nabla_{\!TZ})\wedge{\text{ch}}(\nabla^u_{\!\xi})+{\text{ch}}(\nabla_{\!\eta})
-{\text{ch}}(\nabla_{\!\mathcal H})\]
(where ${\rm{ch}}(\nabla_{\!\eta})={\rm{ch}}(\nabla_{\!\eta^+})-{\rm{ch}}(\nabla_{\!\eta^-})$ and accordingly for ${\rm{ch}}(\nabla_{\!\mathcal H})$). The preceding considerations about functoriality apply here, so
that this transgression form $\int_0^{+\infty}\phi
{\text{Tr}}_s\left(\frac{\partial B_t}{\partial t}
\exp-B^2_t\right)dt$ is functorial by pullbacks on fibered products as in
\eqref{pullbackfunctorial} (if the horizontal subspace of the
source manifold is taken as described at the end of \S\ref{s:geomdata}).
\subsection{The ``families Chern-Simons'' transgression form:}\label{IranzoSarah}
Consider now some vector bundle $\xi$ with connection $\nabla_{\!\xi}$
and hermitian metric $h^\xi$
on $M$, some horizontal tangent vector space $T^H\!M$ and vertical
metric $g^Z$ for the submersion $\pi\colon M\to B$, and some vector
bundles $F^+$ and $F^-$ on $B$ such that (see \cite{MoiPartie1} \S3.1 for a description of $\pi^{\rm{Eu}}_*$)
\[[F^+]-[F^-]=\pi^{\rm{Eu}}_*[\xi]\in K^0_{\text{top}}(B)\]
Put connections $\nabla_{\!F^+}$ on $F^+$ and $\nabla_{\!F^-}$ on $F^-$ and choose
some equivalence class of link $[\ell]$ between $F^+-F^-$ and
some vector bundles $({\mathcal H}^+\oplus\eta^-)-({\mathcal H}^-\oplus\eta^+)$ provided by any suitable data $(\eta^+,\eta^-,\psi)$ (with connections $\nabla_{\!{\mathcal H}}$ and $\nabla_{\!\eta}$). (See \cite{MoiPartie1} \S\S2.1.2 and 2.2.5 for details about links and their Chern-Simons transgression forms).
\begin{definition}\label{deftau}The families Chern-Simons transgression form is
the class modulo exact forms of the following differential form on $B$:
\begin{align*}
\tau(\nabla_{\!\xi},\nabla_{\!TZ},\nabla_{\!F^+},\nabla_{\!F^-},[\ell])=
\int_0^{+\infty}&\phi
{\text{Tr}}_s\left(\frac{\partial B_t}{\partial t}
\exp-B^2_t\right)dt+\\&+\int_Ze(\nabla_{\!TZ})\wedge\widetilde{\rm{ch}}(\nabla_{\!\xi}^u,\nabla_{\!\xi})
+\widetilde{\rm{ch}}([\ell])
\end{align*}
\end{definition}
Of course $\widetilde{\rm{ch}}([\ell])$ is computed with the connections $\nabla_{\!\eta}$, $\nabla_{\!{\mathcal H}}$ and $\nabla_{\!F^\pm}$ (see \cite{MoiPartie1} (29)), and the
form $
\tau(\nabla_{\!\xi},\nabla_{\!TZ},\nabla_{\!F^+},\nabla_{\!F^-},[\ell])$
verifies the following formula (which justifies its name):
\begin{equation}\label{transgressiondesfamillescompleteyeah}
d\tau(\nabla_{\!\xi},\nabla_{\!TZ},\nabla_{\!F^+},\nabla_{\!F^-},[\ell]) =\int_Ze(\nabla_{\!TZ})\wedge{\text{ch}}(\nabla_{\!\xi})+{\text{ch}}(\nabla_{\!F^-})
-{\text{ch}}(\nabla_{\!F^+})
\end{equation}
\begin{lemma*}
$\tau(\nabla_{\!\xi},\nabla_{\!TZ},\nabla_{\!F^+},\nabla_{\!F^-},[\ell])$ does not depend on
$h^\xi$, nor on the function $\chi$ nor on the construction of topological direct image and the 
choice of data used in it, provided the class of link $[\ell]$ is modified
by composition with the canonical link between the obtained representatives of
the topological direct image (as defined in \cite{MoiPartie1} \S3.1.2) when they are changed.
Moreover, $\tau$ is functorial by pullback on fibered products as \eqref{pullbackfunctorial}
provided the horizontal subspaces verify the assumption of the end of \S\eqref{s:geomdata}
\end{lemma*}
$\tau(\nabla_{\!\xi},\nabla_{\!TZ},\nabla_{\!F^+},\nabla_{\!F^-},[\ell])$ of course depends
on the other data in a way which will be precised later.
\begin{proof}
The general principle of the proof will be to use functoriality
to compare restrictions to $B\times\{0\}$ and $B\times \{1\}$ of
some $\widetilde\tau$ constructed as above on a submersion of the form $\pi\times{\rm{Id}}_{[0,1]}\colon
M\times[0,1]\longrightarrow B\times[0,1]$.
The vertical tangent space of $\pi\times{\rm{Id}}_{[0,1]}$ is simply the
pullback to $M\times[0,1]$ of the one of $\pi$, and it will be supposed to be endowed with a pullback metric. Choose some horizontal subspace $T^H\!M$ for $\pi$
and pull it back on $M\times[0,1]$, where it is a suitable horizontal subspace
with respect to $\pi\times{\rm{Id}}_{[0,1]}$. These choices of horizontal subspaces
verify the conditions of the end of \S\ref{s:geomdata} with respect to the
maps $B\times\{0\}\hookrightarrow B\times[0,1]$ and $B\times\{1\}\hookrightarrow B\times[0,1]$ to which we plan to apply functoriality.
Call $\widetilde\nabla_{\!TZ}$ the associated pullback connection on the
vertical tangent bundle of $\pi\times{\rm{Id}}_{[0,1]}$.

Consider some vector bundle $\xi$ on $M$,
with connection
$\nabla_{\!\xi}$. Consider any pair of vector bundles $F^+$ and $F^-$ on $B$
with connections $\nabla_{\!F^+}$ and $\nabla_{\!F^-}$ such that
$[F^+]-[F^-]=\pi^{\rm{Eu}}_*[\xi]\in K^0_{\rm{top}}(B)$, and some equivalence class of link
$[\ell]$ between $F^+-F^-$ and vector bundles obtained from the
families index construction.
Pull back $\xi$ on $M\times[0,1]$ and $F^+$ and $F^-$ on
$B\times[0,1]$ and call $\widetilde\xi$, $\widetilde F^+$ and $\widetilde F^-$
the pullbacks. Call $\widetilde\nabla_{\!\xi}$, $\widetilde\nabla_{\!F^+}$ and
$\widetilde\nabla_{\!F^-}$ the pullback connections on them.

Endow $\widetilde\xi$ with some not necessarily pullback hermitian metric
$\widetilde h^\xi$ and choose any suitable data $(\widetilde\eta^+,\widetilde\eta^-,\widetilde\psi)$
with respect to $\pi\times{\rm{Id}}_{[0,1]}$ to obtain vector bundles
$\widetilde{\mathcal H}^\pm={\rm{Ker}}{\mathcal D}^\pm_{\widetilde\psi}$ on
$B\times[0,1]$ (following the construction of \cite{MoiPartie1} \S3.1.1). One has
\[[\widetilde F^+]-[\widetilde F^-]=
[\widetilde{\mathcal H}^+]-[\widetilde{\mathcal H}^-]=(\pi\times{\rm{Id}}_{[0,1]})^{\rm{Eu}}_*
\widetilde\xi\in K^0_{\rm{top}}(B\times[0,1])\]
$[\ell]$ naturally provides an
equivalence class of link between $F^+-F^-$ and the restrictions to
$B\times\{0\}$ of $(\widetilde{\mathcal H}^+\oplus\eta^-)-(\widetilde{\mathcal H}^-
\oplus\eta^+)$, which
can be extended (by parallel transport along $[0,1]$) to an equivalence class of link $[\widetilde\ell]$ between $\widetilde F^+-\widetilde F^-$ and $(\widetilde{\mathcal H}^+\oplus\widetilde\eta^-)-(\widetilde{\mathcal H}^-\oplus\widetilde\eta^+)$
on the whole $B\times[0,1]$.

Construct the differential form $\widetilde\tau=\tau(\widetilde\nabla_{\!\xi},\widetilde\nabla_{\!TZ},\widetilde\nabla_{\!F^+},\widetilde\nabla_{\!F^-},[\widetilde\ell])$ in the same way as in definition \ref{deftau}
with respect to all these data on $M\times[0,1]$. This must be made using a smooth
function $\widetilde\chi$ on $B\times[0,1)\times{\mathbb R}_+$ vanishing on
$B\times[0,1]\times[0,\varepsilon]$, equal to $1$ on $B\times[0,1]\times[A,
+\infty)$ and increasing with respect to the variable in ${\mathbb R}^+$
as was sketched at the end of \S\ref{s:adapt}. The obtained form $\widetilde\tau$ verifies
\[d\widetilde\tau=\int_Ze(\widetilde\nabla_{\!TZ}){\rm{ch}}(\widetilde\nabla_{\!\xi})
+{\rm{ch}}(\widetilde\nabla_{\!F^-})-{\rm{ch}}(\widetilde\nabla_{\!F^+})\]
where $\int_Z$ stands for integration along the fibres of $\pi\times{\rm{Id}}_{[0,1]}$. Call $\tau_0$ and $\tau_1$ the restrictions of $\widetilde\tau$ to
$B\times\{0\}$ and $B\times\{1\}$ respectively. Integrating this
formula along $[0,1]$ provides that the following differential form on $B$ is exact:
\[d\bigg(\int_{[0,1]}\widetilde\tau\bigg)=\tau_1-\tau_0+\int_{[0,1]}\int_Ze(\widetilde\nabla_{\!TZ}){\rm{ch}}(\widetilde\nabla_{\!\xi})
+\int_{[0,1]}{\rm{ch}}(\widetilde\nabla_{\!F^-})-\int_{[0,1]}{\rm{ch}}(\widetilde\nabla_{\!F^+})\]
but $\widetilde\nabla_{\!TZ}$ and $\widetilde\nabla_{\!\xi}$ are pullback connections on $M\times[0,1]$ for the projection on the second factor $M\times[0,1]\longrightarrow
M$ and accordingly for $\widetilde\nabla_{\!F^+}$ and $\widetilde\nabla_{\!F^-}$ on $B\times[0,1]$,
so that their Chern characters or Euler form are pullback forms, and their
integral along $[0,1]$ vanish. It follows that $\tau_0$ and $\tau_1$ are equal
modulo exact forms.

The $\int_0^{+\infty}\phi
{\text{Tr}}_s\left(\frac{\partial B_t}{\partial t}
\exp-B^2_t\right)dt$ part of $\tau$ is functorial by pullback on fibered products
as in \eqref{pullbackfunctorial} as was remarked at the end of
subsection \ref{s:4.4}. The $\widetilde{\rm{ch}}$ (and $e(\nabla_{\!TZ})$)
parts are also functorial as was remarked just before \cite{MoiPartie1} equation (18) and after \cite{MoiPartie1} (21),
both under the assumption on horizontal subspaces of the end of \S\ref{s:geomdata}.

Thus $\tau_0$ and $\tau_1$ are both regular definitions of
$\tau(\nabla_{\!\xi},\nabla_{\!TZ},\nabla_{\!F^+},\nabla_{\!F^-},[\ell])$
as in definition \ref{deftau}, because the class of link
between $F^+-F^-$ and the restrictions to $B\times\{1\}$
of $(\widetilde{\mathcal H}^+\oplus\widetilde\eta^-)-
(\widetilde{\mathcal H}^-\oplus\widetilde\eta^+)$ is in the
equivalence class of $[\ell]$ (it can be deformed along $[0,1]$
to the one between the restrictions on $B\times\{0\}$).

Thus the independence of the class of $\tau$ modulo exact forms on
$h^\xi$, $\chi$ and the data $\eta^+$, $\eta^-$ and $\psi$ is proved
with the restriction that data $\eta^+$, $\eta^-$ and $\psi$ can be
deformed from one to another (of course hermitian metrics and functions
of the type $\chi$ can always be deformed from one to another).

The last point to check is that
data used to construct the topological direct image
can be deformed from one to another as above (or
almost).

First remark that if $(\eta^+,\eta^-,\psi)$ are suitable data, then
$(\eta^+\oplus\zeta^+,\eta^-\oplus\zeta^-,\psi)$ also are ($\zeta^+$
and $\zeta^-$ are inert excess vector bundles) and give rise to the same
$\tau$. The same is true for $(\eta^+\oplus\zeta,\eta^-\oplus\zeta,
\psi\oplus{\rm{Id}}_\zeta)$ because the extra term appearing in
$\phi
{\text{Tr}}_s\left(\frac{\partial B_t}{\partial t}
\exp-B^2_t\right)$ is the supertrace on $\zeta\oplus\zeta$ of some
${\rm{End}}(\zeta\oplus\zeta)$-valued differential form whose
diagonal terms are equal.

Call ${\mathcal H}^-$ the kernel of ${\mathcal D}^-_\psi$ (with respect
to suitable data $(\eta^+,\eta^-,\psi)$), and $\iota_{{\mathcal H}^-}$
the immersion ${\mathcal H}^-\hookrightarrow({\mathcal E}^-\oplus\eta^-)$.
Put $\widetilde\eta^+=
\eta^+\oplus{\mathcal H}^-$, $\widetilde\eta^-=\eta^-\oplus{\mathcal H}^-$
and $\widetilde\psi=\psi+\cos(\frac\pi2t)\iota_{{\mathcal H}^-}+\sin(\frac\pi2t){\rm{Id}}_{{\mathcal H}^-}$.

On $B\times[0,1]$, the kernel of ${\mathcal D}^+_{\widetilde\psi}$ is constant
and equal to ${\rm{Ker}}{\mathcal D}^+_\psi$, while the kernel of ${\mathcal D}^-
_{\widetilde\psi}$ equals $\big(\sin(\frac\pi2t)\iota_{{\mathcal H}^-}-
\cos(\frac\pi2t){\rm{Id}}_{{\mathcal H}^-}\big){\mathcal H}^-\subset
\big(({\mathcal E}^-\oplus\eta^-)\oplus{\mathcal H}^-\big)$. Applying
the considerations above to this construction proves that the forms $\tau$
constructed from $(\eta^+,\eta^-,\psi)$ (corresponding to $B\times\{0\}$, with
an extra ${\rm{Id}}_{{\mathcal H}^-}$)
and $(\eta^+\oplus{\mathcal H}^-,\eta^-,\psi+\iota_{\mathcal H}^-)$ (
corresponding to $B\times\{1\}$ with an
extra inert copy of ${\mathcal H}^-$ added to $\eta^-$) differ from an exact
form.

Suppose as in the first alinea of \S3.1.2 of \cite{MoiPartie1} that there exists some vector bundle $\lambda$
on $B$ and some bundle map $\varphi\colon\lambda\longrightarrow{\mathcal E}^-\oplus\eta^-$ such that ${\mathcal D}_{\psi+\varphi}^+$ is surjective onto
${\mathcal E}^-\oplus\eta^-$. Consider then $\widetilde\eta^+=\eta^+\oplus\lambda\oplus
{\mathcal H}^-$, $\widetilde\eta^-=\eta^-$ and $\widetilde\psi=\psi+\cos(\frac\pi2t)\varphi+\sin(\frac\pi2t)
\iota_{{\mathcal H}^-}$. Of course ${\mathcal D}^+_{\widetilde\psi}$ is
surjective on $B\times[0,1]$, and its kernel equals $({\rm{Ker}}{\mathcal D}^+_{\psi
+\varphi})\oplus{\mathcal H}^-$ on $M\times\{0\}$ and $({\rm{Ker}}{\mathcal D}^+_\psi)
\oplus\lambda$ on $M\times\{1\}$. Thus applying the above considerations to
this case, proves that $\tau$ constructed from $(\eta^+\oplus{\mathcal H}^-,\eta^-,\psi+\iota_{\mathcal H}^-)$ (corresponding to $M\times\{1\}$
with an inert copy of $\lambda$ added to $\eta^+\oplus{\mathcal H}^-$) and
from $(\eta^+\oplus\lambda,\eta^-,\psi+\varphi)$ (corresponding to $M\times\{0\}$
with an inert copy of ${\mathcal H}^-$ added to $\eta^+\oplus\lambda$)
differ from an exact form; the parallel transport along $[0,1]$
from $({\rm{Ker}}{\mathcal D}^+_{\psi
+\varphi})\oplus{\mathcal H}^-$ to $({\rm{Ker}}{\mathcal D}^+_\psi)
\oplus\lambda$ (following ${\rm{Ker}}({\mathcal D}^+_{\widetilde\psi}\vert_{M\times\{t\}})$) is easily checked to lie in the equivalence class of the link
between $({\rm{Ker}}{\mathcal D}^+_\psi)-{\mathcal H}^-$ and $({\rm{Ker}}{\mathcal D}^+_{\psi+\varphi})-\lambda$ obtained from the exact sequence \cite{MoiPartie1} (40).
Thus the lemma is proved for any change of suitable data of the type of the
first alinea of \S3.1.2 in \cite{MoiPartie1}.

The general case can be checked by applying the above considerations to exactly the same construction on $M\times[0,1]$ as
in the third alinea of \S3.1.2 of \cite{MoiPartie1}. The lemma is thus proved in full generality.
\end{proof}
\subsection{Properties of the Chern-Simons families transgression form:}
\subsubsection{Anomaly formulae:}
The same trick of deformation on $M\times[0,1]$ will be applied further.
One can deform the connection on $\xi$, or the geometry of the fibration, either
the riemannian metric on the fibres or the horizontal tangent subspace $T^H\!M$,
and these deformations can be performed without restriction.
Denote by $\nabla_{\!\xi}^0$ and $\nabla_{\!TZ}^0$ the connections on $\xi$ and
on the vertical tangent subspace corresponding to data on $M\times\{0\}$
and by $\nabla_{\!\xi}^1$ and $\nabla_{\!TZ}^1$ their counterpart on $M\times\{1\}$,
the obtained formula for any couple $(F^+,F^-)$ of bundles on $B$ with
connections $\nabla_{\!F^+}$ and $\nabla_{\!F^-}$ such that $[F^+]-[F^-]=\pi^{\rm{Eu}}_*[\xi]\in
K^0_{\rm{top}}(B)$ is
\begin{equation}\label{anomalie1}\begin{aligned}
\tau&(\nabla_{\!\xi}^1,\nabla^1_{\!TZ},\nabla_{\!F^+},\nabla_{\!F^-},[\ell])-
\tau(\nabla^0_{\!\xi},\nabla^0_{\!TZ},\nabla_{\!F^+},\nabla_{\!F^-},[\ell])
=\\
&\qquad=\int_Z\Big[e(\nabla_{\!TZ}^0)\wedge\widetilde{\text{ch}}(\nabla^0_{\!\xi},
\nabla^1_{\!\xi})
+\widetilde{e}(\nabla_{\!TZ}^0,\nabla^1_{\!TZ})\wedge{\text{ch}}
(\nabla^1_{\!\xi})\Big]
\\&\qquad=\int_Z\Big[\widetilde{e}(\nabla_{\!TZ}^0,\nabla^1_{\!TZ})\wedge{\text{ch}}
(\nabla^0_{\!\xi})+e(\nabla_{\!TZ}^1)\wedge\widetilde{\text{ch}}(\nabla^0_{\!\xi},
\nabla^1_{\!\xi})
\Big]
\end{aligned}\end{equation}
(See for instance \cite{MoiPartie1} (22) for the last equality).

Now one also can change the bundles on $B$ in the following way:
take suitable $(\eta^+,\eta^-,\psi)$ and call ${\mathcal H}^\pm=
{\rm{Ker}}{\mathcal D}^{\nabla_{\!\xi}\pm}_\psi$, endow ${\mathcal H}^+\oplus\eta^-$
and ${\mathcal H}^-\oplus\eta^+$ with any
connections $\nabla^\uparrow$ and $\nabla^\downarrow$. Consider vector bundles $F^+$, $F^-$, $G^+$
and $G^-$ on $B$ such that $[F^+]-[F^-]=[G^+]-[G^-]=\pi^{\rm{Eu}}_*[\xi]
\in K^0_{\rm{top}}(B)$, choose some connections $\nabla_{\!F^+}$, $\nabla_{\!F^-}$,
$\nabla_{\!G^+}$ and $\nabla_{\!G^-}$ on them, and some
links $[\ell_F]$ and $[\ell_G]$ between $F^+-F^-$ or $G^+-G^-$ respectively and
$({\mathcal H}^+\oplus\eta^-)-({\mathcal H}^-\oplus\eta^+)$. Then from the construction of $\tau$
it follows that 
\begin{equation}\label{anomalie2}\begin{aligned}\tau(\nabla_{\!\xi},\nabla_{\!TZ},\nabla_{\!F^+},\nabla_{\!F^-},[\ell_F])&=
\tau(\nabla_{\!\xi},\nabla_{\!TZ},\nabla^\uparrow,\nabla^\downarrow,[{\rm{Id}}])
+\widetilde{\rm{ch}}([\ell_F])\\
&=\tau(\nabla_{\!\xi},\nabla_{\!TZ},\nabla_{\!G^+},\nabla_{\!G^-},[\ell_G])-\widetilde{\rm{ch}}([\ell_G])
+\widetilde{\rm{ch}}([\ell_F])\\
&=\tau(\nabla_{\!\xi},\nabla_{\!TZ},\nabla_{\!G^+},\nabla_{\!G^-},[\ell_G])+\widetilde{\rm{ch}}([\ell_F\circ
\ell_G^{-1}])
\end{aligned}\end{equation}
(see \cite{MoiPartie1} (30)) where of course $\widetilde{\rm{ch}}([\ell_F])$ and $\widetilde{\rm{ch}}([\ell_G])$ are
computed with $\nabla_{\!F^\pm}$ or $\nabla_{\!G^\pm}$ respectively, and $\nabla^\uparrow$ and $\nabla^\downarrow$.

Formulae \eqref{anomalie1} and \eqref{anomalie2} give all the dependence of $\tau$ on its data.
\subsubsection{Partial caracterisation of $\tau$:}
\begin{lemma}\label{lemmederealite}
$\tau(\nabla_{\!\xi},\nabla_{\!TZ},\nabla_{\!F^+},\nabla_{\!F^-},[\ell])$ vanishes if $\nabla_{\!\xi}$ is flat and if $F^+$ and $F^-$ are the sheaf theoretic direct
images of $\xi$ as flat bundles (with associated flat connections as described in \cite{MoiPartie1} \S3.3).

It is additive in the following sense:
let $\xi_1$ and $\xi_2$ be bundles on $M$ with connections $\nabla_{\!\xi_1}$ and $\nabla_{\!\xi_2}$, let $F^+_1$, $F^-_1$, $F^+_2$ and $F^-_2$ be bundles
with connections on $B$
such that $[F^+_1]-[F^-_1]=\pi^{\rm{Eu}}_*[\xi_1]$ and $[F^+_2]-[F^-_2]=\pi^{\rm{Eu}}_*[\xi_2]$
in $K^0_{\text{top}}(B)$. Let $[\ell_1]$ and $[\ell_2]$ be links between
$F^+_1-F^-_1$ and bundles on $B$ obtained from topological direct image construction for $\xi_1$, and correspondingly for $[\ell_2]$. Then, as the topological direct image construction is additive (as direct sum),
the direct sum $\ell_1\oplus\ell_2$ provides an equivalence class
of link between $(F^+_1\oplus F^+_2)-(F^-_1\oplus F^-_2)$ and bundles on
$B$ obtained from topological direct image construction for $\xi_1\oplus\xi_2$, then
\begin{align*}
\tau(\nabla_{\!\xi_1}\oplus\nabla_{\!\xi_2},&\nabla_{\!TZ},\nabla_{\!F^+_1}\oplus\nabla_{\! F^+_2},\nabla_{\!F^-_1}
\oplus \nabla_{\!F^-_2},[\ell_1\oplus\ell_2])=\\
&=\tau(\nabla_{\!\xi_1},\nabla_{\!TZ},\nabla_{\!F^+_1},\nabla_{\!F^-_1},[\ell_1])
+\tau(\nabla_{\!\xi_2},\nabla_{\!TZ},\nabla_{\!F^+_2},\nabla_{\!F^-_2},[\ell_2])
\end{align*}
\end{lemma}
\begin{proof}The vanishing of $\tau$ for flat bundles is a consequence of the first statement of lemma \ref{l:reel}
and of \cite{BismutLott} theorem 3.17:
in the proposed case, the nullity
proved in lemma \ref{l:reel} holds for all $t>0$. This is due to the fact that $B_t=C_t$
when it is computed from the trivial suitable data $(\{0\},\{0\},0)$. In particular, the link $[\ell]$ in the terms $\int_Ze(TZ)\wedge\widetilde{\rm{ch}}(\nabla_{\!\xi},\nabla_{\!\xi}^u)
-\widetilde{\rm{ch}}([\ell])$ is trivial as link, but it links $F^+$ and
$F^-$ endowed with their sheaf theoretic direct image flat connections $\nabla_{\!F^+}$
and $\nabla_{\!F^-}$,
with $F^+$ and $F^-$ endowed with their metric connections $\nabla_{{\!\mathcal H}^+}$ and
$\nabla_{{\!\mathcal H}^-}$ obtained by the
projection on the kernel of the fibral Dirac operator \eqref{limiteinfinie}.
It is proved in \cite{BismutLott} proposition 3.14 that $\nabla_{{\!\mathcal H}^+}=\nabla^u_{\!F^+}$
and accordingly on $F^-$, and in \cite{BismutLott} theorem 3.17 (see also \cite{MoiPartie1} remark 4) that, up to exact forms
\[\widetilde{\rm{ch}}(\nabla_{\!F^+},\nabla_{\!F^+}^u)-\widetilde{\rm{ch}}(\nabla_{\!F^-},\nabla
^u_{\!F^-})=\int_Ze(\nabla_{\!TZ})\widetilde{\rm{ch}}(\nabla_{\!\xi},\nabla_{\!\xi}^u)\]
The nullity of $\tau(\nabla_{\!\xi},\nabla_{\!TZ},\nabla_{\!F^+},
\nabla_{\!F^-},[{\rm{Id}}])$ follows.

The additivity is a direct consequence of the fact that all the construction
of the transgression form is additive for direct sum data, and accordingly
for Chern-Simons transgressions.
\end{proof}
\begin{theorem*}For
bundles $\xi$ with vanishing Chern classes in $H^{\text{even}}(M,
{\mathbb Q})$, the class $\tau$ is the only one which has
the preceding properties of additivity, functoriality by pull-backs on fibered
products, nullity on flat bundles with their sheaf theoretic direct images
and of course the transgression property \eqref{transgressiondesfamillescompleteyeah}.
\end{theorem*}
The anomaly formulae \eqref{anomalie1} and \eqref{anomalie2} are in this case consequences
of these properties (using functoriality and constructions on $B\times[0,1]$).
\begin{proof}
If $\xi$ has vanishing rational Chern classes, then
some finite direct sum $\xi\oplus\xi\oplus\ldots\oplus\xi$ is
topologically trivial on $X$. The anomaly formula then
relates $\tau$ for $\nabla_{\!\xi}\oplus\nabla_{\!\xi}\oplus\ldots
\oplus\nabla_{\!\xi}$
(and any direct sum of copies of direct image representatives)
and for the canonical flat connection on the trivial
bundle with corresponding flat direct image
(for which $\tau$ vanishes because of (iii)). Dividing 
by the number of copies of $\xi$ produces the desired $\tau$.
\end{proof}
\begin{remark*}
One could generalise to bundles $\xi$ whose restrictions to the fibers
of $\pi$ have vanishing rational Chern classes by adding some property
linking $\tau$ for $\xi$ and $\tau$ for
$\xi\otimes\pi^*\zeta$ where $\zeta$ is any bundle on $B$.
To obtain a general caracterisation would need some more: in the
case of double transgression on complex algebraic manifolds, Weng
\cite{Weng} uses the deformation to the normal cone.

One could hope to obtain a caracterisation of $\tau$ modulo the image
of $K^1_{\rm{top}}(B)$ by the Chern character, with no care of links of
bundles on $B$ with someones obtained by analytic families index
construction. However, the fact that one must consider a not controlled finite number of copies of
$\xi$ would prevent to obtain more than a caracterisation modulo rational cohomology.
\end{remark*}
\subsubsection{Reality and Hodge symmetry:}\label{S.I.}
\begin{lemma}\label{oddnull}$\tau(\nabla_{\!\xi},\nabla_{\!TZ},\nabla_{\!F^+},\nabla_{\!F^-},[\ell])$ is a real form if
$\nabla_{\!\xi}$, $\nabla_{\!F^+}$ and $\nabla_{\!F^-}$ respect some hermitian metrics on $\xi$, $F^+$ and $F^-$.

If ${\rm{dim}}Z$ is even, then $\tau(\nabla_{\!\xi}^*,\nabla_{\!TZ},\nabla_{\!F^+}^*,\nabla_{\!F^-}^*,[\ell])
=\overline{\tau(\nabla_{\!\xi},\nabla_{\!TZ},\nabla_{\!F^+},\nabla_{\!F^-},[\ell])}$ (where $\nabla_{\!\xi}^*$, $\nabla_{\!F^+}^*$ and $\nabla_{\!F^-}^*$
are the transpose adjoints of $\nabla_{\!\xi}$, $\nabla_{\!F^+}$
and $\nabla_{\!F^-}$ with respect to any hermitian metrics
on $\xi$, $F^+$ and $F^-$).

If ${\rm{dim}}Z$ is odd, denote by $0$ the connection on the null rank
vector bundle $\{0\}$ on $B$; consider any suitable data $(\eta^+,\eta^-,\psi)$
giving rise to fibral Dirac operator kernels ${\mathcal K}^+$ and ${\mathcal K}^-$, then $\tau(\nabla_{\!\xi},\nabla_{\!TZ},0,0,[\ell_{\mathcal K}^{\{0\}}]^{-1})=0$
(where $[\ell_{\mathcal K}^{\{0\}}]$ is the canonical class of link between
$({\mathcal K}^+\oplus\eta^-)-({\mathcal K}^-\oplus\eta^+)$ and $\{0\}-\{0\}$
obtained in \cite{MoiPartie1} \S3.2.3 and lemma 10).
\end{lemma}
\begin{proof}The reality of $\tau$ for connections which respect some hermitian metrics is a consequence of lemma \ref{l:reel}
and of the fact that Chern-Simons forms are real when computed from connections which respect some hermitian metrics (see \cite{MoiPartie1}
equations (24) and (29) and the comment after (29). If not, one has to use the precised ``reality'' assertion of lemma \ref{l:reel}, and the fact that $\widetilde{\rm{ch}}([\ell])$ is changed into its conjugate when all connections are simultaneously changed into their adjoint transpose).

For the even dimensional case, if $(\eta^+,\eta^-,\psi)$ are suitable
data giving rise to ${\mathcal H}^+$ and ${\mathcal H}^-$, then data
$\big(\eta^+,\eta^-,-(*_Z\oplus{\rm{Id}}_{\eta^-})\circ\psi\circ
(*_Z^{-1}\oplus{\rm{Id}}_{\eta^+})\big)$ give rise to vector bundles
canonically isomorphic to ${\mathcal H}^+$ and ${\mathcal H}^-$
(see \cite{MoiPartie1} lemma 9). It follows from \cite{MoiPartie1} lemma 12, formula \eqref{nablahache} and the definition of $\overline\nabla^u$ in the proof of lemma \ref{pascoole!} that this isomorphism respects the connections
$\nabla_{{\mathcal H}^+}$ and $\nabla_{{\mathcal H}^-}$.

It is proved in \cite{MoiPartie1}
lemma 9 that this isomorphism is compatible with
change of representatives of topological direct image, and the equality
of the lemma follows lemma \ref{l:reel}
(the equation and the reality assertion, or its precised version if the
connections on $\eta^\pm$ do not respect the hermitian metrics $h^\pm$), and the same considerations as just above about reality of Chern-Simons forms (from \cite{MoiPartie1} equations (24), (26), (29)
and the comment after (29), to be adapted as above to a situation where all connections
are simultaneously changed by their adjoint transposes).

In the odd dimensional case, compare the formulae obtained by calculating
$\tau(\nabla_{\!\xi},\nabla_{\!TZ},0,0,[\ell_{\mathcal K}^{\{0\}}]^{-1})$
with some suitable data $(\eta^+,\eta^-,\psi)$ and with associated
``adjoint'' data  $\big(\eta^-, \eta^+,(*_Z\oplus{\rm{Id}}_{\eta^+})\circ
\psi^*\circ(*_Z^{-1}\oplus{\rm{Id}}_{\eta^-})\big)$ (see the end of \S3.2.3 in \cite{MoiPartie1}).
The result follows from the equality between $[\ell_{{\mathcal K}_1}^{\{0\}}]$
and $[-\ell_{{\mathcal K}_0}^{\{0\}}]$ stated at the end of \S3.2.3 in \cite{MoiPartie1}, the fact that connections (on ${\mathcal K}$) are respected by the corresponding isomorphism (from \eqref{nablahache} and lemma \ref{pascoole!} here and lemma 12 in \cite{MoiPartie1}), lemma 10 of \cite{MoiPartie1}, nullity of $e(\nabla_{\!TZ})$, the same remark about transpose adjoint connections and $\widetilde{\rm{ch}}([\ell])$,
the (precised) reality assertion and the second equation in lemma \ref{l:reel} (which alltogether prove that the two calculations of the same object give opposite results).
\end{proof}
\section{Applications:}\label{cinq}
\subsection{Direct image for free multiplicative $K$-theory:}
Take some $(\xi,\nabla_{\!\xi},\alpha)\in \widehat K_{\rm{ch}}(M)$, choose some vector bundles $F^+$ and $F^-$ on $B$ such that $[F^+]-[F^-]=\pi^{\text{Eu}}_*[\xi]\in K^0_{\text{top}}(B)$, endowed with connexions $\nabla_{\!F^+}$, and $\nabla_{\!F^-}$, and put
\begin{equation}\label {E:defimdir}\begin{aligned}
{}&\pi^{\rm{Eu}}_!(\xi,\nabla_{\!\xi},\alpha)=\\&\ =\left(F^+,\nabla_{\!F^+},\int_Ze(\nabla_{\!TZ})\alpha\right)
-\big(F^-,\nabla_{\!F^-},\tau(\nabla_{\!\xi},\nabla_{\!TZ},\nabla_{\!F^+},\nabla_{\!F^-},[\ell])\big)
\end{aligned}\end{equation}
where $[\ell]$ is any equivalence class of link between $F^+-F^-$ and vector bundles produced by families analytic index construction (here $[\ell]$ doesn't matter because only the class of $\tau$ modulo the Chern character of $K^1_{\text{top}}$ will be relevant).
\begin{theorem}\label{thebigone}
The map $\pi^{\rm{Eu}}_!$ is the zero map if ${\rm{dim}}Z$ is odd.
If ${\rm{dim}}Z$ is even, it defines a real morphism from $\widehat K_{\rm{ch}}(M)$ to $\widehat K_{\rm{ch}}(B)$ (with respect to the following ``conjugation involution'': $(E,\nabla,\alpha)\longmapsto(E,\nabla^*,\overline\alpha)$ where $\nabla^*$ is the adjoint transpose to $\nabla$ with respect to any hermitian metric on $E$)

It verifies:

$\widehat{\rm{ch}}\big(\pi_!(\xi,\nabla_{\!\xi},\alpha)\big)=\int_{M/B}e(\nabla_{\!TZ})\widehat{\rm{ch}}(\xi,\nabla_{\!\xi},\alpha)$

${\mathfrak{B}}\big(\pi_!(\xi,\nabla_{\!\xi},\alpha)\big)=\int_{M/B}e(\nabla_{\!TZ})
{\mathfrak{B}}(\xi,\nabla_{\!\xi},\alpha)$

This morphism is compatible with the morphism
\[(E,\nabla_{\!E})\in K^0_{\rm{flat}}\longmapsto(E,\nabla_{\!E},0)\in\widehat K_{\rm{ch}}\]
with the commutative diagram \cite{MoiPartie1} (34), direct images on $K_{\text{top}}^0$, on $K^1_{\text{top}}$, on $K^0_{\text{flat}}$ and the morphism from $\Omega^{\rm{odd}}(M,{\mathbb C})\big/d\Omega^{\rm{even}}(M,{\mathbb C})$
to $\Omega^{\rm{odd}}(B,{\mathbb C})\big/d\Omega^{\rm{even}}(B,{\mathbb C})$
given by $\alpha\mapsto\int_Ze(\nabla_{\!TZ})\alpha$ (integration along the fibre after
product with $e(\nabla_{\!TZ})$).
\end{theorem}
\begin{proof}
\begin{description}

\item[]\underline{Nullity in the odd-dimensional fibre case}:
$\ $ 
This is a direct consequence of the nullity of $e(\nabla_{\!TZ})$
and of the last statement of lemma \ref{oddnull}. All the other properties stated
in the even dimensional fibre case, follow in the odd dimensional case too, but the compatibility with $\pi_!$ on $K^0_{\rm{flat}}$, this last point is proved as in the even dimensional fibre case below.

Thus ${\rm{dim}}Z$
will be supposed even in all the sequel of the proof.
\item[]\underline{Well defining the class of the image}:
$\ $ 
If $G^+$ and $G^-$  with connexions $\nabla_{\!G^+}$ and $\nabla_{\!G^-}$ are as $F^+$ and $F^-$ such that $[G^+]-[G^-]=
\pi^{\text{Eu}}_*[\xi]\in K_0^{\text{top}}(B)$ then from \eqref{anomalie2}
\[\tau(\nabla_{\!\xi},\nabla_{\!TZ},\nabla_{\!F^+},\nabla_{\!F^-},[\ell_F])-\tau(\nabla_{\!\xi},\nabla_{\!TZ},\nabla_{\!G^+},\nabla_{\!G^-},[\ell_G])=\widetilde{\text{ch}}([\ell_F\circ\ell_G^{-1}])
\]
Thus the formula \eqref{E:defimdir} written with $G^+$ and $G^-$ instead of $F^+$ and $F^-$ provides the same class
in $\widehat K_{\rm{ch}}(B)$ (see the relation defining $\widehat K_{\rm{ch}}$ in the introduction just before \eqref{Sarah1}).
\item[]\underline{Independence on the choice of the representative at the source:}
$\ $ 
Suppose that $(\xi,\nabla_{\!\xi},\alpha)=(\xi',\nabla_{\!\xi'},\alpha')\in \widehat K_{\rm{ch}}(M)$, and that $f\colon\xi\to\xi'$ is some $C^\infty$ isomorphism, then
\[\alpha'=\alpha+\widetilde{\text{ch}}(\nabla_{\!\xi},f^*\nabla_{\!\xi'})+\beta\]
where $\beta$ is a closed form lying in the image of $K^1_{\text{top}}(M)$
under the Chern character. Thus if $[F^+]-[F^-]=\pi^{\text{Eu}}_*[\xi]\in K^0_{\text{top}}(B)$
with connexions $\nabla_{\!F^+}$ on $F^+$ and $\nabla_{\!F^-}$
on $F^-$, one has from \eqref{anomalie1} and \eqref{anomalie2}:
\begin{align*}
\tau(\nabla_{\!\xi},\nabla_{\!TZ},\nabla_{\!F^+},&\nabla_{\!F^-},[\ell_\xi])-\tau(\nabla_{\!\xi'},\nabla_{\!TZ},\nabla_{\!F^+},\nabla_{\!F^-},[\ell_{\xi'}])=\\
&=\int_Ze(\nabla_{\!TZ})\wedge\widetilde{\text{ch}}(\nabla_{\!\xi},f^*\nabla_{\!\xi'})
+\widetilde{\text{ch}}(\ell_\xi\circ\ell_{\xi'}^{-1})
\end{align*}
(for any suitable links $[\ell_{\xi}]$ and $[\ell_{\xi'}]$)
so that if for any closed odd degree form $\gamma$, one denotes by $(0,0,\gamma)$
the element of $\widehat K_{\rm{ch}}$ which for any $(E,\nabla,\alpha)\in \widehat K_{\rm{ch}}$ equals $(E,\nabla,\alpha+\gamma)-(E,\nabla,\alpha)$, one obtains
\begin{align*}
\pi^{\text{Eu}}_!&(\xi,\nabla_{\!\xi},\alpha)-\pi^{\text{Eu}}_!(\xi',\nabla_{\!\xi'},\alpha')=\\
&=\left(0,0,\int_Ze(\nabla_{\!TZ})\wedge\Big(\widetilde{\text{ch}}(\nabla_{\!\xi},f^*\nabla_{\!\xi'})+\alpha
-\alpha'\Big)+\widetilde{\text{ch}}(\ell_\xi\circ\ell_{\xi'}^{-1})\right)\\
&=\left(0,0,\int_Ze(\nabla_{\!TZ})\wedge\beta\right)+\left(0,0,\widetilde{\text{ch}}(\ell_\xi\circ\ell_{\xi'}^{-1})\right)
\end{align*}
which vanishes in $\widehat K_{\rm{ch}}(B)$ since $\widetilde{\text{ch}}(\ell_\xi\circ\ell_{\xi'}^{-1})\in{\text{ch}}\big(K^1_{\text{top}}(B)\big)\subset H^{\text{odd}}(B,{\mathbb C})$ and so
does $\int_Ze(\nabla_{\!TZ})\wedge\beta$ because of the cohomological version of Atiyah-Singer families index theorem
for $K^1_{\text{top}}$.

Of course it is straightforward to check that in $\widehat K_{\rm{ch}}(B)$
\[\pi_!^{\text{Eu}}(\xi_1\oplus\xi_2,\nabla_{\!\xi_1}\oplus\nabla_{\!\xi_2},\alpha_1+\alpha_2)
=\pi_!^{\text{Eu}}(\xi_1,\nabla_{\!\xi_1},\alpha_1)+\pi_!^{\text{Eu}}(\xi_2,\nabla_{\!\xi_2},
\alpha_2)\]
from the additivity of $\tau$ for direct sums (last assertion of lemma \ref{lemmederealite}).
\item[]\underline{Reality of the morphism:}
This is a direct consequence of the two first assertions of lemma \ref{oddnull}.
\item[]\underline{Relation concerning $\widehat{\rm{ch}}$:} $\ $ With notations from
\eqref{E:defimdir}
\begin{align*}
\widehat{\rm{ch}}&\big(\pi_*^{\text{Eu}}(\xi,\nabla_{\!\xi},\alpha)\big)=\\&=
{\rm{ch}}(\nabla_{\!F^+})-\int_Ze(\nabla_{\!TZ})\wedge d\alpha-{\rm{ch}}(\nabla_{\!F^-})+d\tau(\nabla_{\!\xi},\nabla_{\!TZ},\nabla_{\!F^+},\nabla_{\!F^-},[\ell])\\
&=\int_Ze(\nabla_{\!TZ})\big({\rm{ch}}(\nabla_{\!\xi})-d\alpha\big)
\end{align*}
because of the transgression relation \eqref{transgressiondesfamillescompleteyeah}.

Note that this relation has the consequence that $\pi_!^{\rm{Eu}}$
sends $MK_0(M)$ to $MK_0(B)$, and
$K^{-1}_{{\mathbb C}/{\mathbb Z}}(M)$ to $K^{-1}_{{\mathbb C}/{\mathbb Z}}(B)$.
\item[]\underline{Relation concerning ${\mathfrak{B}}$:} $\ $
With notations from
\eqref{E:defimdir}
\begin{align*}
{\mathfrak B}&\big(\pi_*^{\text{Eu}}(\xi,\nabla_{\!\xi},\alpha)\big)=\\&=
\widetilde{\rm{ch}}(\nabla_{\!F^+}^*,\nabla_{\!F^+})-\int_Ze(\nabla_{\!TZ})\wedge\alpha+\overline{\int_Ze(\nabla_{\!TZ})\wedge\alpha}
-\widetilde{\rm{ch}}(\nabla_{\!F^-}^*,\nabla_{\!F^-})+\\
&\qquad\qquad\qquad+\tau(\nabla_{\!\xi},\nabla_{\!TZ},\nabla_{\!F^+},\nabla_{\!F^-},[\ell])-
\overline{\tau(\nabla_{\!\xi},\nabla_{\!TZ},\nabla_{\!F^+},\nabla_{\!F^-},[\ell])}
\end{align*}
Of course the connections on $F^+$ and on $F^-$ can be supposed to respect some hermitian metrics on $F^+$ and $F^-$ without changing the formula, and
this makes and vanish the terms $\widetilde{\rm{ch}}(\nabla_{\!F^+}^*,\nabla_{\!F^+})$ and $\widetilde{\rm{ch}}(\nabla_{\!F^-}^*,\nabla_{\!F^-})$.

Now consider any connection $\nabla_{\!\xi}^u$
which respects some hermitian metrics on $\xi$. Then $\tau(\nabla_{\!\xi}^u,\nabla_{\!TZ},\nabla_{\!F^+},\nabla_{\!F^-},[\ell])$ is real whatever $[\ell]$
may be as proved in lemma \ref{lemmederealite}. And from the anomaly
formula \eqref{anomalie1}, one obtains
\begin{align*}
\tau(\nabla_{\!\xi},&\nabla_{\!TZ},\nabla_{\!F^+},\nabla_{\!F^-},[\ell])=\\&=\tau(\nabla_{\!\xi}^u,\nabla_{\!TZ},\nabla_{\!F^+},\nabla_{\!F^-},[\ell])+\int_Ze(\nabla_{\!TZ})\widetilde{\rm{ch}}(\nabla_{\!\xi}^u,\nabla
_{\!\xi})\end{align*}
Now because $\overline{\widetilde{\rm{ch}}(\nabla_{\!\xi}^u,\nabla
_{\!\xi})}=\widetilde{\rm{ch}}(\nabla_{\!\xi}^{u*},\nabla
_{\!\xi}^*)=\widetilde{\rm{ch}}(\nabla_{\!\xi}^u,\nabla
_{\!\xi}^*)$ (see \cite{MoiPartie1} (25))
one gets:
\begin{align*}{\mathfrak B}\big(\pi_!^{\text{Eu}}(\xi,\nabla_{\!\xi},\alpha)\big)&=
-\int_Ze(\nabla_{\!TZ})\wedge\alpha
+{\int_Ze(\nabla_{\!TZ})\wedge\overline\alpha}\, +\\&\qquad+
\int_Ze(\nabla_{\!TZ})\widetilde{\rm{ch}}(\nabla_{\!\xi}^u,\nabla
_{\!\xi})-\int_Ze(\nabla_{\!TZ})\widetilde{\rm{ch}}(\nabla_{\!\xi}^u,\nabla
_{\!\xi}^*)\\
&=\int_Ze(\nabla_{\!TZ})\Big(\widetilde{\rm{ch}}(\nabla_{\!\xi}^*,\nabla_{\!\xi})
-\alpha+\overline\alpha\Big)
\end{align*}
and the relation is proved.

Note that this relation has the consequence that $\pi_!^{\rm{Eu}}$ sends
$\widehat{{\mathbb R}K}_{\rm{ch}}(M)$ to $\widehat{{\mathbb R}K}_{\rm{ch}}(B)$
and $K^{-1}_{{\mathbb R}/{\mathbb Z}}(M)$ to $K^{-1}_{{\mathbb R}/{\mathbb Z}}(B)$.
\item[]\underline{Compatibility with direct images on $K^0_{\text{top}}$, $H^{\text{odd}}$ and
$K^0_{\text{flat}}$:}

The compatibility with the direct image $\pi^{\text{Eu}}_*$ on $K^0_{\text{top}}$ is tautological.

The compatibility with the direct image on odd forms modulo exact forms is
also trivial since $\pi_!^{\text{Eu}}(0,0,\alpha)=\big(0,0,\int_Ze(\nabla_{\!TZ})\wedge\alpha\big)$.

The compatibility with direct image on $K^0_{\text{flat}}$ is a direct consequence of the vanishing of $\tau$ for flat bundles with their sheaf theoretic direct images as stated in lemma \ref{lemmederealite}.
\end{description}
\end{proof}
\subsection{Influence of the vertical metric and the horizontal distribution:}
If geometric data are changed on $M$, namely the vertical riemannian metric $g^Z$ and/or the horizontal subspace $T^H\!M$, this changes the connection $\nabla_{\!TZ}$, and this also changes the morphism $\pi_!^{\rm{Eu}}$.
\begin{lemma}\label{l:anomaliemetrique} Let $\nabla_{\!TZ}$ and $\pi^{\rm{Eu}}_!$
be associated to data $g^Z$ and
$T^H\!M$, let ${g^Z}'$ and ${T^H}'\!M$ be other data and call $\nabla'_{\!TZ}$
and ${\pi_!^{\rm{Eu}}}'$ the associated connection on $TZ$ and morphism from $\widehat K_{\rm{ch}}(M)$
to $\widehat K_{\rm{ch}}(B)$. Then,
for any $(\xi,\nabla_{\!\xi},\alpha)$ one has
\[{\pi^{\rm{Eu}}_!}'(\xi,\nabla_{\!\xi},\alpha)-
\pi_!^{\rm{Eu}}(\xi,\nabla_{\!\xi},\alpha)=\left(0,0,-\int_Z\widetilde e
(\nabla_{\!TZ},\nabla'_{\!TZ})\widehat{\rm{ch}}(\xi,\nabla_{\!\xi},\alpha)\right)
\]
\end{lemma}
\begin{proof} If ${\rm{dim}}Z$ is odd, then both $\pi_!^{\rm{Eu}}$ and ${\pi^{\rm{Eu}}_!}'$ vanish, and $\widetilde e(\nabla_{\!TZ},\nabla'_{\!TZ})$
also does.
If ${\rm{dim}}Z$ is even, it follows from \eqref{anomalie1} that
\[
\tau(\nabla_{\!\xi},\nabla'_{\!TZ},\nabla_{\!F^+},\nabla_{\!F^-},[\ell])-\tau(\nabla_{\!\xi},
\nabla_{\!TZ},\nabla_{\!F^+},\nabla_{\!F^-},[\ell])=\int_Z\!\widetilde e(\nabla_{\!TZ},\nabla'_{\!TZ})\wedge{\text{ch}}(\nabla_{\!\xi})
\]
thus
\begin{align*}
{\pi^{\rm{Eu}}_!}'(\xi,\nabla_{\!\xi},\alpha)&-
\pi_!^{\rm{Eu}}(\xi,\nabla_{\!\xi},\alpha)=\\&=\left(0,0,\int_Z\big(e(\nabla'_{\!TZ})-e(\nabla_{\!TZ})\big)\alpha-\int_Z\widetilde e
(\nabla_{\!TZ},\nabla'_{\!TZ})\wedge{\rm{ch}}(\nabla_{\!\xi})
\right)=\\
&=\left(0,0,\int_Z\widetilde e
(\nabla_{\!TZ},\nabla'_{\!TZ})\wedge\big(-{\rm{ch}}(\nabla_{\!\xi})+d\alpha\big)\right)\end{align*}
this last equality is valid modulo exact forms because
\[d\big(\widetilde e(\nabla_{\!TZ},\nabla'_{\!TZ})\alpha\big)=e(\nabla'_{\!TZ})\alpha
-e(\nabla_{\!TZ})\alpha+(-1)^{{\rm{deg}}\widetilde e(\nabla_{\!TZ},\nabla'_{\!TZ})}\widetilde e(\nabla_{\!TZ},\nabla'_{\!TZ})d\alpha\]
and $\widetilde e(\nabla_{\!TZ},\nabla'_{\!TZ})$ is of degree ${\rm{dim}}Z-1$
which is odd if the fibres are even dimensional.
\end{proof}
If ${\rm{dim}}Z$ is even, and since $\widetilde e(\nabla_{\!TZ},\nabla'_{\!TZ})$ is of degree ${\rm{dim}}Z-1$, it follows that $MK_0$ is the biggest subgroup of $\widehat K_{\rm{ch}}$
on which there is no variation of $\pi_!^{\rm{Eu}}$ when geometric data
$g^Z$ and $T^H\!M$ are changed.
This gives a topological significance to the direct image morphism
$\pi^{\rm{Eu}}_!$ on $MK_0$.
\subsection{Grothendieck-Riemann-Roch for relative $K$-theory:}\label{GRRR}
\begin{theorem}\label{GRRrel}
\[{\mathcal N}_{\text{ch}}\big(\pi_*(E,\nabla_{\!E},F,\nabla_{!F},f)\big)
=\int_Ze(\nabla_{\!TZ})\wedge {\mathcal N}_{\text{ch}}(E,\nabla_{\!E},F,\nabla_{\!F},f)\]
\end{theorem}
\begin{proof}
The anomaly formulae \eqref{anomalie1} and \eqref{anomalie2} yield that
\begin{align*}\tau(\nabla_{\!E},\nabla_{\!TZ},\nabla_{\!\pi_!^+E},\nabla_{\!\pi_!^-E},[{\rm{Id}}])
&-\tau(\nabla_{\!F},\nabla_{\!TZ},\nabla_{\!\pi_!^+F},\nabla_{\!\pi_!^-F},[{\rm{Id}}])=\\
&\qquad=\int_Ze(\nabla_{\!TZ})\widetilde{\rm{ch}}(\nabla_{\!E},f^*\nabla_{\!F})
-\widetilde{\rm{ch}}([\ell_{\pi_*f}])
\end{align*}
(see \cite{MoiPartie1} definition 15 for the definition of $[\ell_{\pi_*f}]$).
 
Both $\tau$ vanish (see lemma \ref{lemmederealite}), and that the
right hand side vanishes is exactly the desired result in view of the
definition of ${\mathcal N}_{\rm{ch}}$ (just before corollary \ref{CorollSarah}). (Note that both $\widetilde{\rm{ch}}$
are closed forms so that only the cohomology class of $e(\nabla_{\!TZ})$ is
needed here).
\end{proof}
\begin{remark*}
In this statement, (as in the caracterisation of $\tau$ and contrarily to the result on $\widehat K_{\rm{ch}}$), the links between representatives of the topological direct image are to be
considered.
\end{remark*}
\subsection{Some comments about a theorem of J.-M. Bismut on CCS classes:}
Let $(E,\nabla_{\!E})$ be some flat vector bundle on $M$. Its class in
$K^0_{\rm{top}}(M)$ minus its rank, is torsion, so that there exists some (nonunique) integer
$k$ such that $kE$ (which means $E\oplus E\oplus\ldots\oplus E$ with $k$ copies of $E$) is a topologically trivial vector bundle on $M$. Let $f\colon{\mathbb C}^{k{\rm{rk}}E}\overset\sim\longrightarrow kE$ be some trivialisation,
then the form $\frac1k\widetilde{\rm{ch}}\big(d,f^*(k\nabla_{\!E})\big)$ is closed and its cohomology class defines some class (known as the ``Cheeger-Chern-Simons class'') $CCS(E,\nabla_{\!E})\in H^{\rm{odd}}
(M,{\mathbb C})/H^{\rm{odd}}(M,{\mathbb Q})$. The quotient is taken to obtain independency on $k$ and on the trivialisation $f$ chosen.
(Here $d$ is the canonical trivial connection on the trivial vector bundle ${\mathbb C}^{k{\rm{rk}}E}$ and $k\nabla_{\!E}$ the canonical direct sum connection on $kE$).

Of course $CCS(E,\nabla_{\!E})$ is the class of $\frac1k{\mathcal N}_{\rm{ch}}({\mathbb C}^{k{\rm{rk}}E},d,kE,k\nabla_{\!E},f)$ modulo $H^{\rm{odd}}(M,{\mathbb Q})$.
\begin{theorem*}(Bismut \cite{BismutChine})
\[CCS(\pi_!^+E,\nabla_{\!\pi_!^+E})-CCS(\pi_!^-E,\nabla_{\!\pi_!^-E})
=\int_Ze(\nabla_{\!TZ})\, CCS(E,\nabla_{\!E})\]
\end{theorem*}
The ``imaginary part'' of this theorem first appeared in \cite{BismutLott} theorem 3.17, it was reproved by X. Ma and W. Zhang
in \cite{MaZhang} and is contained in the result on the
class ${\mathfrak B}$ of theorem \ref{thebigone}, since
$i{\mathfrak{Im}}CCS(E,\nabla_{\!E})=\frac12{\mathfrak B}(E,\nabla_{\!E},0)$ as can be checked from the additivity of
${\mathfrak B}$ for direct sums and from \cite{MoiPartie1} (25) and commentary between (24) and it. However the proof here (of the compatibility of
$\pi_!$ on $K^0_{\rm{flat}}$ and $\pi_!^{\rm{Eu}}$ on $\widehat K_{\rm{ch}}$)
heavily relies on theorem 3.17 in \cite{BismutLott} (see lemma \ref{lemmederealite} above).

A particular case of this theorem is that if ${\mathbb C}$ stands for the
rank one trivial vector bundle on $M$ with its canonical trivial connection $d$, then
\[CCS(\pi_!^+{\mathbb C},\nabla_{\!\pi_!^+{\mathbb C}})
-CCS(\pi_!^-{\mathbb C},
\nabla_{\!\pi_!^-{\mathbb C}})=0\]
which is not trivial since $\pi_!^+{\mathbb C}$
and $\pi_!^-{\mathbb C}$ need not be trivial flat bundles on $B$.

However, once this is granted, the full result will follow from
theorem \ref{GRRrel} since the class modulo $H^{\rm{odd}}(B,{\mathbb Q})$ of $\frac1k{\mathcal N}_{\rm{ch}}\big(\pi_*({\mathbb C}^{k{\rm{rk}}E},d,kE,k\nabla_E,f)\big)$ will be equal to
\[CCS(\pi_!^+E,\nabla_{\!\pi_!^+E})-CCS(\pi_!^-E,\nabla_{\!\pi_!^-E})-CCS(\pi_!^+{\mathbb C},\nabla_{\!\pi_!^+{\mathbb C}})
+CCS(\pi_!^-{\mathbb C},
\nabla_{\!\pi_!^-{\mathbb C}})\]

\end{document}